\documentclass[a4paper,12pt]{article}
\usepackage[latin1]{inputenc}
\usepackage{a4,amsmath,amssymb,amscd,array,booktabs}
\usepackage[all]{xy}

\newtheorem{prop}{Proposition}[section]
\newtheorem{thm}[prop]{Theorem}
\newtheorem{rem}[prop]{Remark}
\newtheorem{lema}[prop]{Lemma}
\newtheorem{coro}[prop]{Corollary}

\newcommand{\Q}{\mathbb Q}
\newcommand{\Qbar}{\overline{\Q}}

\newcommand{\R}{\mathbb R}
\newcommand{\Z}{\mathbb Z}
\newcommand{\N}{\mathbb N}
\newcommand{\F}{\mathbb F}
\newcommand{\PP}{\mathbb P}
\newcommand{\Fbar}{\overline{\F}}
\newcommand{\C}{\mathbb C}
\newcommand{\Gal}{\mathrm{Gal}}
\newcommand{\G}{\mathrm{G}}

\newcommand{\HH}{\mathbb H}
\newcommand{\GL}{\mathrm{GL}}
\newcommand{\PGL}{\mathrm{PGL}}
\newcommand{\SL}{\mathrm{SL}}
\newcommand{\PSL}{\mathrm{PSL}}

\newcommand{\eps}{\varepsilon}
\newcommand{\To}{\longrightarrow}
\newcommand{\Hom}{\operatorname{Hom}}

\newcommand{\Aut}{\operatorname{Aut}}

\newcommand{\qed}{\hfill $\square$}

\title{A moduli approach to quadratic $\Q$-curves
realizing projective mod~$p$\, Galois representations}
\author{Julio Fern\'andez}

\begin{document}
\maketitle

% ABSTRACT ------------------------------------------------------------
\begin{abstract}
For a fixed odd prime $p$\, and a representation $\varrho$\, of the absolute
Galois group of $\Q$ into the projective group $\PGL_2(\F_p)$, we provide
the twisted modular curves whose rational points supply the quadratic
$\Q$-curves of degree $N$ prime to $p$ that realize $\varrho$ through the 
Galois action on their \mbox{$p$-torsion} modules. The modular curve to twist
is either the fiber product of the modular curves $X_0(N)$ and $X(p)$ or
a certain quotient of Atkin-Lehner type, depending on the value of $N$
mod $p$. For our purposes, a special care must be taken in fixing
rational models for these modular curves and in studying their automorphisms.
By performing some genus computations, we obtain from Faltings' theorem
some finiteness results on the number of quadratic $\Q$-curves of a given
degree $N$ realizing $\varrho$.
\end{abstract}

%-----------------------------------------------------------------------
\section{\large Introduction}\label{pgl}

Throughout we fix an odd prime $p$\, and an algebraic closure $\Qbar$ of $\Q$.
For a subfield~$L$ of $\Qbar$, we denote by~$\G_L$ the absolute Galois group $\Gal(\Qbar/L)$. 
We say that a non-CM elliptic curve defined over a quadratic field $k$ is a 
\emph{$\Q$-curve of degree~$N$} if there is a cyclic isogeny of degree~$N$ from the curve 
to its Galois conjugate. The $p$-torsion of a $\Q$-curve~$E_{/k}$ of degree $N$ prime to $p$\, 
gives rise to a Galois representation
$$\varrho_E\,\colon\,\G_\Q\,\To\,\PGL_2(\F_p)$$
whose conjugacy class is an invariant of the isomorphism class of $E$ and whose restriction to~$\G_k$ is 
the projective representation obtained from the usual Galois action on the
\mbox{$p$-torsion} points of the curve. The procedure to obtain $\varrho_E$ is 
detailed in Section~\ref{Weil}. The determinant of $\varrho_E$ 
draws off two different cases, which we call \emph{cyclotomic} and \emph{non-cyclotomic},
and which correspond to $N$\, being a square mod~$p$ or not, respectively.
These two cases rule most of the structure and contents of the rest of sections.

The situation just described raises the following inverse problem: find the \mbox{$\Q$-curves}
of degree $N$ realizing a given projective mod $p$\, Galois representation
$$\varrho\,\colon\,\G_\Q\,\To\,\PGL_2(\F_p)\mathrm{,}$$
namely those whose $p$-torsion points give rise to $\varrho$\, in the above sense. The aim of 
this paper is to explain in detail how to produce the moduli spaces whose rational points yield 
the solutions to this problem. Henceforth the term \emph{rational} stands for \emph{$\Q$-rational}.

The moduli spaces that we provide are either twists of the modular curve $X(N,p)$ obtained as 
the fiber product of the curves $X_0(N)$ and~$X(p)$, in the non-cyclotomic case, or twists 
of a certain Atkin-Lehner quotient $X^+(N,p)$ in the cyclotomic case.
In Section~\ref{complex} we analyze the structure of the subgroup $\mathcal{W}(N,p)$ of
automorphisms on $X(N,p)$ extending the group generated by the
Atkin-Lehner involution~$w_{\scriptscriptstyle N}$ on~$X_0(N)$.
In Section~\ref{Q} we fix a suitable rational model for
$X(N,p)$ and then describe the Galois action on $\mathcal{W}(N,p)$. In
order to do this, we need first to study the action of $\mathcal{W}(N,p)$ on the
non-cuspidal points of the curve. The last two sections explain how to
obtain the twisted curves whose non-cuspidal \mbox{non-CM} rational points give the $\Q$-curves
of degree $N$ realizing $\varrho$. Section \ref{cyclotomic}
is devoted to the cyclotomic case and Section \ref{non-cyclotomic} to the non-cyclotomic case. 
They also include some finiteness results obtained from Faltings' theorem and from some
genus computations performed in Section \ref{complex}.

Our moduli approach turns out to be quite effective in the non-cyclotomic case, since the 
quadratic field of definition for the possible $\Q$-curves realizing the 
representation $\varrho$ is uniquely determined and, for every fixed degree~$N$\!, 
we just need one twist~$X(N,p)_\varrho$\,. For an explicit application we refer to 
\cite{FGL}, where a plane quartic model is provided for the genus-three case $X(5,3)_\varrho$\,.

In the cyclotomic case, one should instead consider two twists 
$X^+(N,p)_\varrho$\,, $X^+(N,p)'_\varrho$\, whose rational points
include the cyclic isogenies of degree $N$ between elliptic curves over $\Q$ realizing 
$\varrho$. One may also approach the problem by adding a given
quadratic field $k$\, as extra data: the $\Q$-curves of degree $N$ 
defined over~$k$ realizing $\varrho$ are given by the non-cuspidal non-CM rational points on 
two other twisted curves $X(N,p)_{\varrho,\,k}$\,, $X(N,p)'_{\varrho,\,k}$\,.

%------------------------------------------------------------------------
\section{\large Projective mod $p$\, Galois representations
realized by $p$-admissible $\Q$-curves}\label{Weil}

The aim of this section is to review the construction of the 
representation
$$\varrho_E\,\colon\,\G_\Q\,\To\,\PGL_2(\F_p)$$
attached to a ($p$-admissible) $\Q$-curve $E$ and to compute its
determinant. Up to some minor points, the section is mostly a
reformulation of known facts that can mainly be found in
\cite{El-Sk} and go back to \cite{Rib92}. The particular case of
quadratic $\Q$-curves is written down in \cite{Se-Da}
using the ideas of \cite{Sh78}.

Let $E$\, be a \emph{$\Q$-curve}. By this we mean a non-CM elliptic 
curve defined over a number field $L$ and with an isogeny
$$\lambda_\sigma\,\colon\,^\sigma\! E\,\To\,E$$
for every $\sigma$ in $\G_\Q$. Without loss of generality, we
always take $\lambda_\sigma$ equal to~$\lambda_\tau$ whenever
$\sigma$~and~$\tau$ restrict to the same embedding of~$L$ 
into~$\Qbar$, and one might also assume the isogenies $\lambda_\sigma$
to be cyclic. We suppose here that the $\Q$-curve $E$ is
\emph{$p$-admissible}, namely that the isogenies~$\lambda_\sigma$
can be chosen so that $p$\, does not divide the degree of
any of them.

For an isogeny $\varphi\,\colon E'\!\To\!E$, let us write $\varphi^{-1}$ for
the element $(\deg\varphi)^{-1}\!\otimes\widehat\varphi$ 
in \mbox{$\Q\otimes\Hom(E,E')$}, where $\widehat{\varphi}$\, is the dual
isogeny of~$\varphi$. 
Since $E$ has no~CM, any isogeny $E'\!\To\! E$ differs from $\varphi$\, by a 
rational number. Thus, the \mbox{$2$-cocycle} of $\G_\Q$
$$c_E\,\colon\,(\sigma,\tau)\,\longmapsto\,\lambda_\sigma\,
{^\sigma\!\lambda_\tau}\,\lambda_{\sigma\tau}^{-1}$$ 
takes values in $\Q^*$\!.
Let $\alpha$ be a \emph{splitting map} for the $2$-cocycle $c_E$ viewed 
inside the trivial cohomology group $H^2(\G_\Q,\Qbar^{\,*})$, that is, a 
continuous map $\G_\Q\!\To\!\Qbar^{\,*}$ satisfying
$$\lambda_\sigma\,{^\sigma\!\lambda_\tau}\,\lambda_{\sigma\tau}^{-1} \, = \,
\alpha(\sigma)\,\alpha(\tau)\,\alpha(\sigma\tau)^{-1}$$
for all $\sigma,\tau$ in $\G_\Q$. By taking degrees, one deduces that
the map \,\,$\sigma\mapsto\alpha(\sigma)^2\!/\!\deg\lambda_\sigma$ is a Galois 
character. In particular, the values taken by $\alpha$ are algebraic integers 
prime to~$p$. So there exist a finite extension $\F_\alpha$ of $\F_p$ and 
a mod $p$\, reduction map $\widetilde\alpha\,\colon\G_\Q\!\To\!\F_\alpha^{\,*}$\,
obtained from a fixed embedding of $\Qbar$ into a fixed algebraic closure~$\Qbar_p$ 
of $\Q_p$.

Consider now the $\F_\alpha$\!--linear action of $\G_\Q$ on 
$\F_\alpha\!\otimes_{\F_p}\!E[p]$ given by
$$\left(\sigma,\,1\!\otimes\!P\right)\,\longmapsto\,
\widetilde\alpha(\sigma)^{-1}\!\otimes\!
\lambda_\sigma({^\sigma\! P})\mathrm{.}$$ 
By means of the choice of a basis for the $\F_p$\!--module $E[p]$, this action 
produces a linear representation
$$\rho_{E,\,\alpha}\,\colon\,\G_\Q\,\To\,\F_\alpha^{\,*}\,\GL_2(\F_p)$$
defined up to conjugation by matrices in $\GL_2(\F_p)$. The
corresponding projective Galois representation \,$\varrho_E$\, is
actually given by the induced action
$$\left(\sigma,\,C\right)\,\longmapsto\,\lambda_\sigma(^\sigma\! C)$$
on the projective line
$$\PP\left(\,E[p]\,\right) \, = \, 
\left\{\,C\subset E[p] \ \big| \ C\simeq\F_p\,\right\}\mathrm{.}$$ 
This projective representation \,$\varrho_E$\, depends on neither the $p$-admissible system 
of isogenies~$\lambda_\sigma$ nor the splitting map $\alpha$. Further,
the following proposition shows that \,$\varrho_E$\, is an invariant of the 
\emph{$p$-admissible isogeny class} of~$E$.

\begin{prop}\label{isogeny}
Let $E'$ be an elliptic curve over $\Qbar$ and \,$\varphi\,\colon E'\!\To\!E$ be
an isogeny of degree prime to $p$. Then \,$\varrho_{E'}=\varrho_E$.
\end{prop}

\noindent\emph{Proof.}\, Let $\widehat\varphi$\, be the dual isogeny of 
$\varphi$\, and let $\lambda_\sigma$ and $\alpha$ be as before. Consider 
the $2$-cocycle $c_{E'}$ attached to
the $p$-admissible system of isogenies
$\widehat\varphi\,\,\lambda_\sigma\, {^\sigma\!}\varphi$\,
for the $\Q$-curve $E'$\!. Then $\alpha\deg\varphi$\, is a splitting map 
for $c_{E'}$ whose reduction mod~$p$
takes values in the same finite field \,$\F$ as $\widetilde\alpha$. 
The isomorphism $E'[p]\!\To\!E[p]$ induced by~$\varphi$\, extends naturally to an 
isomorphism $\F\!\otimes E'[p]\!\To\!\F\!\otimes E[p]$ that is compatible with 
the corresponding $\F$--linear actions of $\G_\Q$. So $\rho_{E,\,\alpha}$ and 
$\rho_{E'\!,\,\alpha\deg\varphi}$ are conjugated by a matrix in 
$\GL_2(\F_p)$, and the result follows. \qed

\begin{rem}{\rm \,
The (conjugacy class of the) representation $\rho_{E,\,\alpha}$ is
the linear mod~$p$\, representation obtained from the Galois action on the abelian 
variety of $\GL_2$-type attached in \cite{Rib92} to the $\Q$-curve $E$ and
the splitting map~$\alpha$.
Moreover, any lifting of $\varrho_E$ into $\GL_2(\Fbar_p)$ is of the form 
$\rho_{E,\,\alpha}$ for some splitting map \,$\alpha$\, for $c_E$.}
%indeed, a continuous map $\G_\Q\to\Qbar^{\,*}$ is a splitting map for $c_E$ 
%if and only if it is the twist of $\alpha$ by a Galois character
\end{rem}

Note that the restriction of $\varrho_E$ to $\G_L$ is the 
projective representation
$$\overline\rho_E\,\colon\,\G_L\,\To\,\PGL_2(\F_p)$$
obtained from the usual Galois action on the $p$-torsion points of~$E$.
In terms of number fields, this provides 
the fixed field of $\varrho_E$ with the following property: its
composite with $L$ is the splitting field of the modular polynomial 
$\Phi_p(j_E\,;\,X)$ over~$L$, where $j_E$ stands for the $j$-invariant
of the elliptic curve $E$. Whenever $L$~is normal over $\Q$ and 
\,$\overline\rho_E$\, is surjective, this property singles out the 
fixed field of~$\varrho_E$ among all Galois extensions of~$\Q$ with 
group~$\PGL_2(\F_p)$. 

We recall that the determinant of \,$\overline\rho_E$ is the 
restriction to $\G_L$ of the quadratic Galois character
$$\eps\,\colon\,\G_\Q\,\To\,\F_p^{\,*}/{\F_p^{\,*}}^2\,\simeq\,\{\pm 1\}$$ 
obtained from the mod $p$\, cyclotomic character $\chi$. The fixed 
field of \,$\eps$\, is the only quadratic field $k_p\!=\!\Q(\sqrt{\pm p\,}\,)$ 
inside the $p$-th cyclotomic extension of $\Q$. Let us now show that the
projective representation \,$\varrho_E$\,
is odd by first computing the determinant of a lifting.

\begin{prop}\label{det-linear}
The determinant of \,$\rho_{E,\,\alpha}$ is the product of the mod $p$ cyclotomic
character $\chi$ and the character \,$\G_\Q\!\To\!\F_\alpha^{\,*}$\, defined by
\,\,$\sigma\mapsto\,\deg\lambda_\sigma/\,\widetilde\alpha(\sigma)^{2}$\!.
\end{prop}
\noindent\emph{Proof.}\, By virtue of the properties of the Weil pairing
${\langle\,\cdot\,, \cdot\,\rangle}_{E,\,p}$\,, the equalities
$${\langle\,\lambda_\sigma(^\sigma\! P),\,\lambda_\sigma(^\sigma\!
Q)\,\rangle}_{E,\,p} \ = \ {\langle\,^\sigma\!
P,\,\widehat\lambda_\sigma\lambda_\sigma(^\sigma\!
Q)\,\rangle}_{^\sigma\! E,\,p} \ = \ \left(\,{\langle\, P,\,
Q\,\rangle}_{E,\,p}\,\right)^{\chi(\sigma)\deg\lambda_\sigma}$$
hold for any two points $P, Q$ in $E[p]$ and any $\sigma$ in
$\G_\Q$. Whenever $[P,Q]$ is a basis of~$E[p]$, so is 
$\left[\lambda_\sigma(^\sigma\! P),\,\lambda_\sigma(^\sigma\! Q)\right]$.
Moreover, the left-hand term in the above equalities is the power of 
${\langle\, P,\, Q\,\rangle}_{E,\,p}$ to the determinant of the basis 
change. Therefore, the result follows from the definition 
of $\rho_{E,\,\alpha}$\,.\qed

\begin{coro}\label{odd}
The Galois representation $\varrho_E$ is odd.
\end{coro}

\noindent\emph{Proof.}\, Since $\varrho_E$ does not depend on the
$p$-admissible system of isogenies $\lambda_\sigma$ chosen, we can take 
$\lambda_\sigma$ as the identity for all $\sigma$ in $\G_L$. Fix a
splitting map $\alpha$ for the $2$-cocycle $c_E$ obtained
from the isogenies $\lambda_\sigma$. Note that the restriction of
$\alpha$ to~$\G_L$ is then a Galois character.
Write $\varsigma$ for the complex conjugation in $\G_\Q$
obtained by fixing an embedding $\Qbar\!\To\!\C$. 
We must prove the equality 
$\det\rho_{E,\,\alpha}(\varsigma)=-1$. Consider the isogeny 
\,$\lambda_\varsigma\,{^\varsigma\lambda_\varsigma}\,\colon E\!\To\!E$. 
It is given, on the one hand, by multiplication by $\alpha(\varsigma)^2$ 
and, on the other hand, by multiplication by $\pm\deg\lambda_\varsigma$. 
All we have to do is to pin down the latter sign. As a complex elliptic 
curve, $E$\, is isomorphic to $E_z\!=\C/(\Z + z\,\Z)$\, for some $z$ in the
complex upper-half plane $\HH$.  Through this isomorphism, 
$\lambda_\varsigma\,{^\varsigma\lambda_\varsigma}$ translates into the 
isogeny $E_z\!\To\!E_z$ induced by multiplication by 
$\delta\,{^\varsigma}\delta$ for some $\delta$ in $\C^*$\!. 
So the above sign is positive and thus
$\widetilde\alpha(\varsigma)^2=\deg\lambda_\varsigma$ in~$\F_p$. 
Since $\chi(\varsigma)=-1$, the result follows from 
Proposition \ref{det-linear}.\qed\\

Let \,$deg\,\colon\G_\Q\!\To\!\Q^*\!/{\Q^*}^2$\, be the degree character
induced by any $p$-admissible system of isogenies 
$\lambda_\sigma\colon{^\sigma\!}E\!\To\!E$. Then, consider the 
\mbox{mod $p$}\, degree character 
$$deg_p\,\colon\,\G_\Q\,\To\,\F_p^{\,*}/{\F_p^{\,*}}^2
\,\simeq\,\{\pm 1\}$$ 
obtained from \,$deg$\, by composition with the natural map 
\ $\Q^*\!/{\Q^*}^2\!\To\!\Q_p^*/{\Q_p^*}^2$\!. 
The following statement is a 
straightforward consequence of Proposition \ref{det-linear}.

\begin{coro}\label{det-proj}
The determinant of \,$\varrho_E$ is the product \,$\eps\,deg_p$.
\end{coro}
%Notice that $\det\varrho_E$ is trivial only when $k=k_p$ and 
%$N$ is not a square mod $p$.

\begin{rem}\label{deg}{\rm \ 
If the map \,$deg$\, is not trivial, its fixed field \,$K_{deg}$\, is a composite 
of quadratic fields \,$\Q(\sqrt{a_1}\,), \dots, \Q(\sqrt{a_m}\,)$,
where $2^m$ is the degree of \,$K_{deg}$ over $\Q$.
For every \,$l=1,\dots,m$\,, take \,$\sigma_l$\, in~$\G_\Q$ restricting to 
the non-trivial automorphism of $K_{deg}$ that fixes 
$\sqrt{a_h}$\, for $h\not=l$. 
Then, the map \,$deg_p$ is the product of the quadratic Galois characters attached to 
the extensions $\Q(\sqrt{a_l}\,)$ for which 
$\deg\lambda_{\sigma_l}$ is not a square mod $p$.
%As shown in \cite{Elk}, there is always a $\Q$-curve over $K_{deg}$ isogenous to $E$ 
%and for which the isogenies to its Galois conjugates can be taken to have kernel of 
%squarefree order hnece cyclic.
}
\end{rem}

We say that a projective mod $p$\, Galois representation
$$\varrho\,\colon\,\G_\Q\,\To\,\PGL_2(\F_p)$$
is \emph{realized by} a ($p$-admissible) $\Q$-curve $E$ if $\varrho_E = \varrho$, 
where this equality makes only sense up to conjugation in $\PGL_2(\F_p)$. 
The rest of sections are devoted to the particular case of quadratic $\Q$-curves. 
Assume $\varrho$\, to be realized by a $p$-admissible $\Q$-curve of degree $N$\!, 
that is, by a non-CM elliptic curve defined over a quadratic field and with a cyclic 
isogeny to its Galois conjugate of degree $N$ prime to $p$. From Corollary~\ref{det-proj} 
and Remark~\ref{deg}, $\varrho$ has determinant $\eps$ if and only if $N$ is a square 
mod $p$, and otherwise any $\Q$-curves of degree~$N$ realizing~$\varrho$ must be defined 
over the fixed field of the quadratic character \,$\eps\det\varrho$. We refer to the 
first case ($N$ square mod~$p$) as the \emph{cyclotomic case}, and to the second one 
($N$~non-square mod~$p$) as the \emph{non-cyclotomic case}. 

%------------------------------------------------------------------------
\section{\large Automorphisms of the modular curve~$X(N,p)$}\label{complex}

Let $N>1$ be an integer prime to $p$. Let $X_0(N)$, $X(p)$ and $X(1)$ be the 
modular curves attached to the congruence subgroups $\Gamma_0(N)$, $\Gamma(p)$ 
and $\SL_2(\Z)$, respectively. We denote by $X(N,p)$ the modular curve attached 
to the congruence subgroup \,$\Gamma_0(N)\cap\Gamma(p)$, namely the fiber product 
of $X_0(N)$ and $X(p)$ over $X(1)$\,:
$$
\xymatrix{
 X(N,p) \ar@{->}[rd] \ar@{->}[d] &   \\
 X_0(N) \ar@{->} [rd] & X(p) \ar@{->}[d] \\
 & X(1)}$$
The aim of this section is to introduce a certain group \,$\mathcal{W}(N,p)$\,
of automorphisms on~$X(N,p)$. We also compute the genus of this curve.

As a complex curve, $X(N,p)$ is a Galois covering of $X_0(N)$ with group 
\,$\mathcal{G}(N,p)$ given by the quotient $\Gamma_0(N)/\!\pm\Gamma_0(N)\cap\Gamma(p)$.
Since the mod $p$\, reduction map 
$\SL_2(\Z)\!\To\!\SL_2(\F_p)$ induces the exact sequence
$$1 \ \To \ \pm\Gamma_0(N)\cap\Gamma(p) \ \To \ \Gamma_0(N) \ \To
\ \PSL_2(\F_p) \ \To \ 1\mathrm{,}$$
there is a canonical isomorphism
$$\mathcal{G}(N,p)\ \simeq \ \PSL_2(\F_p)\mathrm{.}$$
We recall that \,${\mathcal G}(N,p)$\, consists of the automorphisms \,$g$\, on $X(N,p)$ 
for which the following diagram commutes:
$$\xymatrix{
 X(N,p) \ar@{->}[rd] \ar@{->}[rr]^g & & X(N,p) \ar@{->}[ld]  \\
 & X_0(N) & }$$

Let $w_{\scriptscriptstyle N}$ be the Atkin-Lehner involution on $X_0(N)$ and denote 
by $X^+(N)$ the corresponding quotient. For any integers $a, b, c, d$ satisfying 
\mbox{$a\,d\,N - b\,c\,p^2 = 1$} and~\,$d\equiv\pm 1\, (\mathrm{mod} \ p)$, the action of the
matrix
$$\left(\begin{matrix}a\,N & b\,p\\ c\,p\,N & d\,N \end{matrix}\right)$$
on the complex upper-half plane $\HH$ defines an automorphism \,$\vartheta$\, on $X(N,p)$ 
\emph{extending}~$w_{\scriptscriptstyle N}$, namely making the following diagram commutative:
$$\xymatrix{
 X(N,p) \ar@{->}[r]^\vartheta \ar@{->}[d] & X(N,p) \ar@{->}[d] \\
 X_0(N) \ar@{->}[r]^{{w}_N} & X_0(N)}$$
Indeed, one can check that the above matrix lies in the normalizer of 
\,$\Gamma_0(N)\cap\Gamma(p)$ inside $\PSL_2(\R)$. Hence, the covering $X(N,p)\!\To\!X^+(N)$ has as many
automorphisms as its degree, which means that it is a Galois covering. Let~$\mathcal{W}(N,p)$ 
denote its automorphism group:
$$\xymatrix{ X(N,p) \ar@{->}[d]\ar@/^1pc/@{.}[1,0]^{\mathcal{G}(N,\,p)}
 \ar@/_2pc/@{.}[2,0]_{\mathcal{W}(N,\,p)}\\
X_0(N) \ar@{->}[d] \\
X^+(N)
}$$
The group \,$\mathcal{W}(N,p)$\, contains \,$\mathcal{G}(N,p)$\, as a subgroup of index two whose 
complement consists of the automorphisms on $X(N,p)$ extending $w_{\scriptscriptstyle N}$.

\begin{prop}\label{W}
The group $\mathcal{G}(N,p)$ is a direct factor of \,$\mathcal{W}(N,p)$ if and only if~\,$N$ 
is a square mod $p$. More precisely, the structure of \,$\mathcal{W}(N,p)$ is as follows:
\begin{itemize}
\item[$\cdot$] In the cyclotomic case, there is a unique involution \,$w$ on $X(N,p)$ such that
$$\mathcal{W}(N,p) \ = \ \mathcal{G}(N,p)\,\times\,\langle w\rangle\mathrm{.}$$
\item[$\cdot$] In the non-cyclotomic case,
$$\mathcal{W}(N,p) \ \simeq \ \PGL_2(\F_p)\mathrm{.}$$
\end{itemize}
In the first case, the quotient curve $X(N,p)/w$ is a Galois covering of $X^+(N)$ with group 
$\mathcal{G}(N,p)$. In the second case, the quotient of $X(N,p)$ by an involution in
$\mathcal{W}(N,p)$ is never a Galois covering of $X^+(N)$.
\end{prop}

\noindent\emph{Proof.}\, Viewed as the quotient $\SL_2(\F_p)/\{\pm 1\}$, the group
$\PSL_2(\F_p)$ is generated by the matrices
$$T\,=\,\left(\begin{matrix}1 & 1 \\ 0 & 1\end{matrix}\right)\!\mathrm{,}
\ \ \ \ \ \ \ \ \
U\,=\,\left(\begin{matrix}1 & 0 \\
1 & 1\end{matrix}\right)\!\mathrm{.}$$
%which have order $p$\, and satisfy the relation \,$T\,U^{p-2}\,T = U^2$\!.
On the other hand, the determinant \,$\GL_2(\F_p)\!\To\!\F_p^{\,*}$\, induces an exact sequence
$$1 \ \To \ \PSL_2(\F_p) \ \To \ \PGL_2(\F_p) \ \stackrel{\det}\To
\ \F_p^{\,*}/{\F_p^{\,*}}^2 \ \To \ 1\mathrm{,}$$
so that $\PSL_2(\F_p)$ can be identified with a subgroup of $\PGL_2(\F_p)$ of index two whose 
complementary subgroups are those generated by a conjugate of the matrix
$$V\,=\,\left(\begin{matrix}0 & -v \\ 1 & \,\,\,0\end{matrix}\right)\!\mathrm{,}$$
where $v$ is a non-square in $\F_p^{\,*}$. Since one has the relations $V\,T = U^{-v^{-1}}\,V$ 
and $V\,U = T^{-v}\,V$\!, a system of generators for $\PGL_2(\F_p)$ is given by $V$ and either $T$ or~$U$\!.
Now, $\mathcal{G}(N,p)$ is generated by the automorphisms defined by the matrices
$$T_{\scriptscriptstyle N}\,=\,\left(\begin{matrix}1 & 1 \\ 0 & 1
\end{matrix}\right)\!, \ \ \ \ \ \ \ \ \
U_{\scriptscriptstyle N}\,=\,\left(\begin{matrix}1 & 0 \\
\tilde N\,N & 1\end{matrix}\right)$$
in $\Gamma_0(N)$, where \,$\tilde N\in\Z$\, is any inverse of $N$ mod $p$.
To give a complementary subgroup for \,$\mathcal{G}(N,p)$\, inside
$\mathcal{W}(N,p)$, let us consider separately the two possibilities
for $N$ mod $p$\,:
\begin{itemize}

\item[$\cdot$] If $N$ is a square mod $p$, then it is also a square
mod $p^2$. Let \,$a, b$\, be any integers satisfying \,$a^2 N - b\,p^2 = 1$.
Then, the matrix
$$Z_{\scriptscriptstyle N}\,=\,
\left(\begin{matrix}a\,N & b\,p\\ p\,N & a\,N\end{matrix}\right)$$
defines an involution $w$ on $X(N,p)$ extending
$w_{\scriptscriptstyle N}$. Moreover, $w$ commutes with the
automorphisms defined by $T_{\scriptscriptstyle N}$ and
$U_{\scriptscriptstyle N}$, so it generates a direct cofactor
of \,$\mathcal{G}(N,p)$ inside $\mathcal{W}(N,p)$. The uniqueness
of $w$ comes from the fact that $\PSL_2(\F_p)$ has trivial center.

\item[$\cdot$] If $N$ is not a square mod $p$, then neither is $\tilde N$\!.
Moreover, the matrix
$$V_{\scriptscriptstyle N}\,=\,\left(\begin{matrix}0 & -1\\N & \,\,\,0
\end{matrix}\right)\!\mathrm{,}$$
which defines an involution on $X(N,p)$ extending $w_{\scriptscriptstyle N}$,
satisfies the relations
\,$V_{\scriptscriptstyle N}\,T_{\scriptscriptstyle N}\,=
\,U_{\scriptscriptstyle N}^{\,-N}\,V_{\scriptscriptstyle N}$\, and
\,$V_{\scriptscriptstyle N}\,U_{\scriptscriptstyle N}\,=
\,T_{\scriptscriptstyle N}^{\,-\tilde N}\,V_{\scriptscriptstyle N}$\, inside
$\mathcal{W}(N,p)$. These are precisely the relations that the matrix $V$\!, 
for \,$v$\, equal to~$\tilde N$ mod~$p$, satisfies with the generators \,$T, U$ of 
$\PSL_2(\F_p)$. Hence, the group $\mathcal{W}(N,p)$ is isomorphic to $\PGL_2(\F_p)$.
\end{itemize}
The last assertion in the statement follows from the group structure of~\,$\mathcal{W}(N,p)$\,: 
in the first case, the subgroup $\langle w\rangle$ is normal, while in the second case
$\mathcal{W}(N,p)$ has no normal subgroups of order two because it has trivial center. \qed

\begin{rem}\label{canonical}{\rm \ 
The matrices $Z_{\scriptscriptstyle N}$ and
$V_{\scriptscriptstyle N}$ in the proof of Proposition \ref{W} have
determinant $N$\!. Thus, in the same way as the automorphisms
in $\mathcal{G}(N,p)$ are defined by matrices in $\Gamma_0(N)$ acting
on $\HH$, the
automorphisms on $X(N,p)$ extending $w_{\scriptscriptstyle N}$ are
defined by matrices in $\mathrm{M}_2(\Z)$ with determinant $N$ and hence
lying in $\GL_2(\F_p)$ when reduced mod $p$.
%Any two such matrices defining the same automorphism
%in~$\mathcal{W}(N,p)$ differ by a matrix in $\pm\Gamma_0(N)\cap\Gamma(p)$.
So we have a mod $p$\, reduction map
$$\mathcal{W}(N,p)\,\To\,\PGL_2(\F_p)$$
whose restriction to $\mathcal{G}(N,p)$ is the canonical isomorphism
onto $\PSL_2(\F_p)$. In the non-cyclotomic case, this map is
the isomorphism $\mathcal{W}(N,p)\simeq\PGL_2(\F_p)$
constructed in the proof of Proposition \ref{W}. We keep this \emph{canonical}\, isomorphism
throughout the rest of the paper.}
\end{rem}

\begin{rem}\label{involutions}{\rm \ 
In the non-cyclotomic case, all involutions on $X(N,p)$
extending~$w_{\scriptscriptstyle N}$ are conjugated inside
$\mathcal{W}(N,p)$. Hence, their defining matrices
in $\mathrm{M}_2(\Z)$ can be obtained conjugating the
matrix~$V_{\scriptscriptstyle N}$ in the proof of Proposition \ref{W}
by matrices in $\Gamma_0(N)$. So they can be chosen to be of the
form
$$\left(\begin{matrix}a\,N & \,b\\ c\,N & -a\,N \end{matrix}\right)
\!\mathrm{,}$$
where \,$a, b, c$\, are integers satisfying \,$a^2N + b\,c = -1$.
This fact is used in the proof of Proposition \ref{w-action}.
}
\end{rem}

In the cyclotomic case, let us write $X^+(N,p)$ for the quotient of $X(N,p)$ 
by the only involution $w$ in the center of the group
$\mathcal{W}(N,p)$. To conclude this section, we give a formula for the
genus of $X(N,p)$ and compute the values of~$N$ and $p$\, for which the
curves $X(N,p)$ and $X^+(N,p)$ have genus zero or one. In the proof of
Proposition \ref{genus}, we recall the description of the
cusps of $X_0(N)$. We refer to \cite{Go} for this, as well
as for the action of the Atkin-Lehner involutions on the set of cusps.
Both things are used in the proof of Proposition \ref{X+(N,p)rational}.

\begin{prop}\label{genus}
The genus of the modular curve $X(N,p)$ is
$$1 \ + \ \frac{\,\psi(N)\,p\,(p^2-1)\,}{24} \ - \
\frac{\,p^2-1\,}{4}\!\sum_{0<n\,|N}
\!\varphi\left(\,\left(n,N/n\right)\,\right)\mathrm{,}$$
where $(a,b)$, $\varphi(r)$ and $\psi(N)$ are the usual notations for
the greatest common divisor of the integers $a$ and $b$, the
order of the group $(\Z/r\,\Z)^*$ and the index of\,~$\Gamma_0(N)/\{\pm 1\}$
in $\PSL_2(\Z)$, respectively.
\end{prop}

\noindent\emph{Proof.}\, The number of cusps of $X_0(N)$ is
$\sum\varphi(h_n)$, where the sum is taken over the positive divisors~$n$
of~$N$\!, and $h_n$ stands for $(n,N/n)$. For every divisor $n$, there is
exactly one cusp for each integer $m$ in a system of representatives in $\Z$
of the group $(\Z/h_n\,\Z)^*$\!. We just take the integer $m=1$ whenever
$\varphi(h_n)=1$. Any such integer~$m$ can be chosen prime to $n$,
and the corresponding cusp is then represented by the rational number $m/n$.
The ramification degree of this cusp over $X(1)$ is~$N/(n\,h_n)$.
On the other hand, the cusps of $X(p)$ have ramification degree~$p$ over $X(1)$. 
Thus, since $N$ is prime to~$p$, the cusps of $X(N,p)$
also have ramification degree~$p$ over $X_0(N)$. Moreover, $X(p)$ has no
elliptic points, so neither has $X(N,p)$. Lastly, the degrees of the
coverings $X(N,p)\!\To\!X_0(N)$ and $X_0(N)\!\To\!X(1)$ are \,$p\,(p^2-1)/2$\,
and~$\psi(N)$, respectively. Hence, the proposition follows from the
Hurwitz formula applied to the map $X(N,p)\!\To\!X(1)$. \qed

%\begin{rem}{\rm If $N$ is the product of $l$ different primes
%\,$p_1,\dots, p_l$\,, then the formula for the genus of $X(N,p)$ becomes
%$$1 \ + \ \frac{\,p^2-1\,}{4}\,
%\left(\frac{\,p\,(p_1+1)\,\cdots\,(p_l+1)\,}{6}
%\ - \ 2^l \right)\mathrm{.}$$
%}
%\end{rem}

\begin{coro}\label{X(N,p)rational-elliptic}
The modular curve $X(N,p)$ has genus greater than one, except for
the genus-zero case $X(2,3)$ and the elliptic case $X(4,3)$.
\end{coro}

\noindent\emph{Proof.}\, Since the genera of $X(p)$ and $X_0(N)$
are greater than one for $p>5$ and $N>49$, respectively, one only
has to check the values that Proposition \ref{genus} yields in the
remaining cases. \qed

\begin{lema}\label{X+(pN)rational}
For an odd prime $p$ and an integer $N>1$ prime to $p$, consider the
Atkin-Lehner involution $w_{\scriptscriptstyle N}$ on the modular
curve $X_0(pN)$. The only pairs $(N,p)$ for which the quotient curve
$X_0(pN)/w_{\scriptscriptstyle N}$ has genus zero are $(2,3)$, $(4,3)$,
$(5,3)$, $(8,3)$, $(11,3)$, $(2,5)$, $(4,5)$ and $(3,7)$.
\end{lema}

\noindent\emph{Proof.}\, For every integer $D>71$, the modular curve
$X_0(D)$ has positive genus and is neither elliptic nor hyperelliptic \cite{Ogg}.
For each odd prime $p$\, and each integer $N$ prime to $p$\, such that $pN\leq 71$, 
one can then use the formulae in \cite{Kl}
or the tables \cite{STNB} to conclude the lemma. \qed

\begin{prop}\label{X+(N,p)rational}
The curve $X^+(N,p)$ has genus greater than one, except for the genus-zero
case $X^+(4,3)$.
\end{prop}

\noindent\emph{Proof.}\, The involution $w$, which is defined by the
matrix $Z_{\scriptscriptstyle N}$ in the proof of Proposition \ref{W},
restricts to the Atkin-Lehner involution $w_{\scriptscriptstyle N}$
on $X_0(pN)$, so it induces a Galois covering 
$X^+(N,p)\!\To\!X_0(pN)/w_{\scriptscriptstyle N}$.
On the other hand, the cusps of~$X_0(pN)$ that ramify on $X(N,p)$ are 
those of the form \,$m/n$\, with $p$\, dividing~$n$, and the ramification
degree is always $p$\, (cf. the proof of Proposition~\ref{genus}). 
In particular, the Hurwitz formula implies
that $X_0(pN)/w_{\scriptscriptstyle N}$ has genus zero whenever $X^+(N,p)$
has genus less than two. By Lemma~\ref{X+(pN)rational}, 
the only pairs $(N, p)$, with~$N$ prime to $p$\,
and square mod $p$, for which $X_0(pN)/w_{\scriptscriptstyle N}$ has
genus zero are $(4,3)$ and $(4,5)$. In the first case, the involution $w$
fixes the cusp $1/2$, so $X^+(4,3)$ is a genus-zero quotient of the
elliptic curve $X(4,3)$. Let us now study the second case, for which we 
consider the following commutative diagram:
$$\xymatrix{
 X(4,5) \ar@{->}[rd]^2 \ar@{->}[d]^{10} &   \\
 X_0(20) \ar@{->}[rd]^2 & X^+(4,5) \ar@{->}[d]^{10} \\
 & X_0(20)/w_{\scriptscriptstyle 4}
}$$
The only ramified points of the covering $X(4,5)\!\To\!X_0(20)$ are cusps. 
Moreover, it can be checked that the points lying above the two cusps \,$1/2,1/10$\,
fixed by $w_4$ are also fixed by the involution $w$. Thus, the only ramified cusps 
of the covering \,$X^+(4,5)\!\To\!X_0(20)/w_4$\, are the points above $1/5$ and~$1/10$, 
all of them with ramification degree $5$. Then, the Hurwitz formula shows that 
there must be ten more ramified points, necessarily with  ramification 
degree $2$\, and lying above the two non-cuspidal points on $X_0(20)$ fixed by $w_4$,
hence the genus of~$X^+(4,5)$ is four.
Notice that there are no other ramified points because the number of points on $X_0(20)$
fixed by $w_4$ is exactly four (cf. \cite{Kl} or \cite{STNB}).\qed

%---------------------------------------------------------------------
\section{\large A rational model for the modular curve~$X(N,p)$}\label{Q}

This section deals with the rationality over $\Q$ for the curve $X(N,p)$ as well 
as for the automorphism group $\mathcal{W}(N,p)$ introduced in the previous section:
we fix a certain rational model for $X(N,p)$ that makes the automorphisms in 
$\mathcal{W}(N,p)$ be defined over $k_p$. Recall that $k_p$ stands for the only 
quadratic field inside the $p$-th cyclotomic extension of~$\Q$. We denote by $\zeta_p$
the root of unity $e^{2\pi i/p}$.

Since $X(N,p)$ is the fiber product of the modular curves $X(p)$ and $X_0(N)$ over~$X(1)$, 
a rational model for the first curve is determined by fixing rational models for the other 
three curves. Recall that the function field of $X(1)$ is generated over $\Q$ by the 
elliptic modular function $j$. For $X_0(N)$, consider the canonical rational model given by 
the function field \,$\Q(j,j_{\scriptscriptstyle N})$, where $j_{\scriptscriptstyle N}$ is the 
modular function defined by $j_{\scriptscriptstyle N}(z)\!=\!j(Nz)$ for $z$ in the complex
upper-half plane~$\HH$. As for~$X(p)$, 
the rational model that we fix satisfies the following property: its extension to $k_p$\, gives by 
specialization over an elliptic curve $E$\, in $X(1)(\Qbar)$ the fixed field of the projective 
mod $p$\, Galois representation $\overline\rho_E$ attached to the \mbox{$p$-torsion} points of $E$. 
This model for $X(p)$ is obtained as the following particular case of a general procedure that 
follows Section II.3 in \cite{Li} and Section~2 in~\cite{Ma77}.

Fix a non-square $v$ in $\F_p^{\,*}$ and take a matrix $V$ in $\GL_2(\F_p)$ of order two in
$\PGL_2(\F_p)$ and with $\det(V)=v$. Without risk of confusion, we often identify the matrix $V$\!, 
up to a sign, with its image in $\PGL_2(\F_p)$. Define $H_V$ as the inverse image in $\GL_2(\F_p)$ 
of the subgroup generated by $V$ in $\PGL_2(\F_p)$\,:
$$H_V\,=\,\F_p^{\,*}\,\cup\,\F_p^{\,*}\,V\!\mathrm{.}$$
Up to conjugation, $H_V$ is the only subgroup of $\GL_2(\F_p)$ containing the center~$\F_p^{\,*}$ 
and reducing inside $\PGL_2(\F_p)$ to a complementary subgroup of $\PSL_2(\F_p)$. The group $H_V$
defines, as shown in the following diagram, a rational model $X_V(p)$ whose 
$\Q$-isomorphism class does not depend on the choice of the matrix $V$:
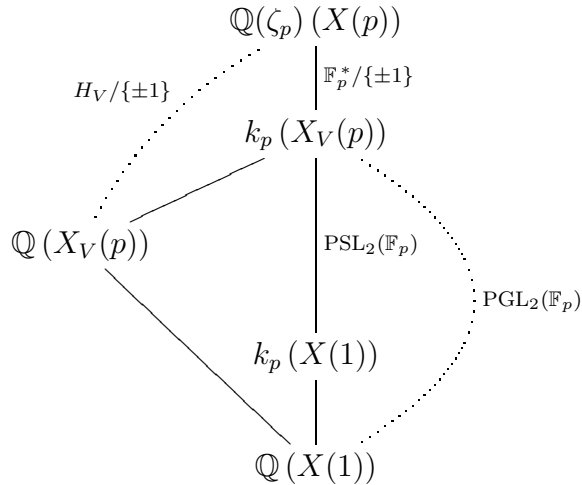
$$\xymatrix{
    & \Q(\zeta_p)\left(X(p)\right)
    \ar@{-}[d]^{{\F_p^{\,*}}/\{\pm 1\}}
    \ar@/_1pc/@{.}[2,-1]_{H_V/\{\pm 1\}}  \\
    & k_p\left(X_V(p)\right)
    \ar@{-}[ld] \ar@{-}[dd]^{\PSL_2(\F_p)}
    \ar@/^5pc/@{.}[3,0]^{\PGL_2(\F_p)}\\
    \Q\left(X_V(p)\right) \ar@{-}[rdd] & \\
    & k_p\left(X(1)\right) \ar@{-}[d] \\
    & \Q\left(X(1)\right)
}$$
Here \,$\Q(\zeta_p)\left(X(p)\right)$\, stands for the field of modular functions for $\Gamma(p)$ whose 
Fourier expansions have coefficients in $\Q(\zeta_p)$. This is a Galois extension of \,$\Q\left(X(1)\right)$ 
with group $\GL_2(\F_p)/\{\pm 1\}$ and the function field of $X_V(p)$ is then the fixed field by the
subgroup $H_V/\{\pm 1\}$. 
%Note that $k_p\left(X_V(p)\right)$ is a Galois extension of $\Q\left(X(1)\right)$ with group $\PGL_2(\F_p)$.

Although we should denote by~$X_V(N,p)$ the rational model for $X(N,p)$ obtained 
from~$X_V(p)$, we just write $X(N,p)$ for simplicity. Without loss of generality, 
we always take the above non-square $v$\, equal to~$N^{-1}$ mod $p$\, in the non-cyclotomic case.  
Note that the map $X(N,p)\!\To\!X_0(N)$ is defined over~$\Q$ and that the function
field $k_p\left(X(N,p)\right)$ is a Galois extension of $\Q\left(X_0(N)\right)$ with 
group~$\PGL_2(\F_p)$. In particular, the Galois action on the automorphism 
group $\mathcal{G}(N,p)$\, factors through $\Gal(k_p/\Q)$.

The non-cuspidal complex points on $X(N,p)$ are in bijection with the isomorphism 
classes of triples
$$\left(E,\,C,\,[T_1,T_2]_V\right)\!\mathrm{,}$$
where $E$ is a complex elliptic curve, $C$ is a cyclic subgroup of $E(\C)$ of order~$N$\!, 
$[T_1,T_2]$ is a basis for $E[p]$ and $[T_1,T_2]_V$ is the corresponding orbit inside 
$E[p]\times E[p]$ by the action of $H_V$. Here $H_V$ is viewed as a subgroup of automorphisms 
of $E[p]$ through the isomorphism \mbox{$\GL_2(\F_p)\simeq\Aut(E[p])$} fixed by the basis 
$[T_1,T_2]$, so that
$$[T_1,T_2]_V \,=\, \left\{\,[r\,T_1, r\,T_2],\,[r\,T_1, r\,T_2]\,V \
\big| \ r\in\F_p^{\,*}\,\right\}\mathrm{.}$$
Two triples of the form $\left(E,\,C,\,[T_1,T_2]_V\right)$ are isomorphic if
there is an isomorphism between the corresponding elliptic curves interchanging the
cyclic subgroups and the $H_V$--orbits.
%sending one cyclic subgroup to the other and one $H_V$--orbit to the other.

This bijection is compatible with the usual Galois actions. Thus, a point on~$X(N,p)$\, given
by a triple as above with \,$j_E\not= 0, 1728$\, is defined over a number field~$L$\, if and only if 
the elliptic curve~$E$ is defined over~$L$, the subgroup $C$ is $\G_L$-invariant and the image
of the linear Galois representation
$$\rho_E\,\colon\,\G_L\,\To\,\GL_2(\F_p)$$
attached to $E[p]$ lies inside a conjugate of the subgroup~$H_V$.

We can always assume that the basis $[T_1,T_2]$ in a triple $\left(E,\,C,\,[T_1,T_2]_V\right)$ 
is, inside the corresponding $H_V$--orbit, the only one up to a sign that is sent to $\zeta_p$
by the Weil pairing. The Galois action on the non-cuspidal points of $X(N,p)$ 
should then be written accordingly: an automorphism \,$\sigma$\, of \,$\C$ takes any
such a triple to that given by the elliptic curve \,${^\sigma\!}E$, the subgroup 
\,${^\sigma\!}C$\, and the $H_V$--orbit of either the basis 
$[r^{-1}\,\,{^\sigma}T_1, r^{-1}\,\,{^\sigma}T_2]$ or the basis 
$[(v\,r)^{-1}\,\,{^\sigma}T_1, (v\,r)^{-1}\,\,{^\sigma}T_2]\,V$\!, depending on whether  
\,$^\sigma\!\zeta_p=\zeta_p^{\,r^2}$\, or \,$^\sigma\!\zeta_p=\zeta_p^{\,v\,r^2}$\,
for some $r$ in $\F_p^{\,*}$, respectively.

The action of the automorphism group $\mathcal{G}(N,p)$ on the non-cuspidal points 
of~$X(N,p)$, and then the Galois action on $\mathcal{G}(N,p)$, are stated in
Proposition~\ref{G(N,p)-action} and Corollary \ref{gal-on-G(N,p)}, respectively. 
The symbol \,$\hat{}$\, stands henceforth for the matrix (anti)involution given by
$$\hat{M} = \left(\begin{matrix}0 & 1\\ 1 & 0\end{matrix}\right)
{^\mathrm{t}\!}M \left(\begin{matrix}0 & 1\\ 1 & 0\end{matrix}\right)
\!\mathrm{,}$$
where ${^\mathrm{t}\!}M$ is the transpose of the matrix $M$.
Alternatively, it can be defined as follows:
$$M = \left(\begin{matrix}a & b\\ c & d\\\end{matrix}\right)
\ \longmapsto \ \hat{M}=\left(\begin{matrix}d & c\\ b & a\end{matrix}\right)\!\mathrm{.}$$

\begin{prop}\label{G(N,p)-action}
An automorphism in \,$\mathcal{G}(N,p)$\, represented through the canonical isomorphism 
\,$\mathcal{G}(N,p)\simeq\PSL_2(\F_p)$\, by a matrix $\gamma$ in~$\SL_2(\F_p)$ takes a point
$\left(E,\,C,\,[T_1,T_2]_V\right)$ on $X(N,p)$ to the point given
by the elliptic curve~$E$, the subgroup $C$ and the \mbox{$H_V$--orbit} of
the $p$-torsion basis $[T_1,T_2]\,\hat{\gamma}$.
\end{prop}

\noindent\emph{Proof.}\, Take any matrix $\left({a\atop c}\,{b\atop d}\right)$ in $\SL_2(\Z)$
reducing mod $p$ to $\gamma$. The triple $\left(E,\,C,\,[T_1,T_2]_V\right)$ is isomorphic 
to one of the form
$$\left(E_z\,,\,\langle 1/N\rangle,\,\left[\,1/p,\,z/p\,\right]_V\right)$$
for some $z$ in $\HH$, where $E_z$ stands for the complex elliptic curve defined by the lattice
$\Z + z\,\Z$. The automorphism in the statement sends the pair given by $z$ to that given by
$$z' \ = \ \frac{\,a\,z + b\,}{\,c\,z + d\,}\,\mathrm{.}$$
Then, the endomorphism of $\C$ defined by multiplication by $c\,z + d$
extends to an isomorphism $E_{z'}\!\To\!E_z$ that preserves the subgroup $\langle 1/N\rangle$
and sends the basis $[\,1/p,\,z'/p\,]$ of $E_{z'}[p]$ to the basis
$$[\,(d + c\,z)/p,\,(b + a\,z)/p\,] \, = \, [\,1/p,\,z/p\,]\,\hat{\gamma}$$
of $E_z[p]$, so the result follows. \qed

\begin{coro}\label{gal-on-G(N,p)}
An automorphism in \,$\mathcal{G}(N,p)$\, represented through the canonical isomorphism
\,$\mathcal{G}(N,p)\simeq\PSL_2(\F_p)$\, by a matrix $\gamma$ in $\SL_2(\F_p)$ is sent by 
the non-trivial element in $\Gal(k_p/\Q)$ to the automorphism in \,$\mathcal{G}(N,p)$ 
corresponding to the matrix \,$\hat{V}\gamma\,\hat{V}$ in $\PSL_2(\F_p)$.
\end{coro}

\noindent\emph{Proof.}\, Denote by $g$ the automorphism represented by the matrix $\gamma$. 
Take any element $\sigma$ in~$\G_\Q$ such that $^\sigma\!\zeta_p=\zeta_p^{\,v}$ and 
let $\gamma_{\sigma}$ be a matrix in $\SL_2(\F_p)$ representing the automorphism
${^\sigma\!}g\in\mathcal{G}(N,p)$. Take also any point~$P$ in~$X(N,p)(\Qbar)$ given by a 
triple $\left(E,\,C,\,[T_1,T_2]_V\right)$ 
with~\mbox{$j_E\not=0, 1728$}. The identity
\mbox{\,${^\sigma\!}\left(\,g\,(P)\,\right) = {^\sigma\!}g\,(^\sigma\! P)$}
leads, by means of Proposition~\ref{G(N,p)-action}, to an
automorphism of the elliptic curve~$^{\sigma\!}E$\, interchanging the
$H_V$--orbits of the \mbox{$p$-torsion} bases
\,$[v^{-1}\,\,{^\sigma}T_1, v^{-1}\,\,{^\sigma}T_2]\,\hat{\gamma}\,V$\, and
\,$[v^{-1}\,\,{^\sigma}T_1, v^{-1}\,\,{^\sigma}T_2]\,V\hat{\gamma}_{\sigma}$.
This yields the identity \,$\hat{\gamma}_{\sigma}=V\hat{\gamma}\,V$\, 
in~$\PSL_2(\F_p)$. \qed\\

From now on, we fix as follows an involution \,$w$\, on $X(N,p)$ 
extending the Atkin-Lehner involution $w_{\scriptscriptstyle N}$ on $X_0(N)$.
Recall that $\mathcal{W}(N,p)$ stands for the group of the Galois covering
$X(N,p)\!\To\!X^+(N)$, where $X^+(N)$ is the quotient of~$X_0(N)$ by
$w_{\scriptscriptstyle N}$. In the cyclotomic case, we take as $w$ the
only involution in the center of~$\mathcal{W}(N,p)$\, (cf. Proposition \ref{W})
and denote by $\sqrt{N}$ a square root of $N$ mod~$p$.
In the non-cyclotomic case, we take as $w$ the involution corresponding
to the matrix~$\hat V$ through the canonical isomorphism
\,$\mathcal{W}(N,p)\simeq\PGL_2(\F_p)$\, (cf. Remark~\ref{canonical}).
Recall that, in the second case, $\det(V)$ is taken to be $N^{-1}$ mod $p$.

\begin{prop}\label{w-action}
The involution \,$w$\, takes a point
\,$\left(E,\,C,\,[T_1,T_2]_V\right)$\, on \,$X(N,p)$ to the point given
by the elliptic curve $E/C$, the subgroup $E[N]/C$ and the \mbox{$H_V$--orbit}
of the image in $E/C$ of the following $p$-torsion basis:
\begin{itemize}
\item[$\cdot$] \,$[\sqrt{N}^{\,-1}\,T_1, \sqrt{N}^{\,-1}\,T_2]$\, in the cyclotomic case.
\item[$\cdot$] \,$[T_1, T_2]\,V$\, in the non-cyclotomic case.
\end{itemize}
\end{prop}

\noindent\emph{Proof.}\, According to Remark \ref{involutions} and the
proof of Proposition \ref{W}, the involution~$w$ is always defined by
the action on $\HH$ of a matrix in $\mathrm{M}_2(\Z)$ of the form
$$\left(\begin{matrix}a\,N & b\\ c\,N & d\,N \end{matrix}\right)\!\mathrm{,}$$
with $a\,d\,N - b\,c = 1$. Denote by \,$\gamma$\, the reduction mod $p$\,
of this matrix. We now proceed as in the proof of Proposition \ref{G(N,p)-action}\,: 
the given triple is isomorphic to one of the form
$$\left(E_z\,,\,\langle 1/N \rangle,\,\left[\,1/p,\,z/p\,\right]_V\right)$$
for some $z$ in $\HH$, where $E_z$ stands for the complex elliptic curve defined by the lattice
$\Z + z\,\Z$.
% $E_z=\C/(\Z + z\,\Z)$.
The involution $w$ sends the triple given by $z$ to that given by
$$z' \ = \ \frac{\,a\,N\,z + b\,}{\,c\,N\,z + d\,N\,}\,\mathrm{.}$$
Then, the endomorphism of $\C$ defined by multiplication by \,$c\,z + d$\,
extends to an isomorphism
$$E_{z'}\ \To\ \ E_z\big/\langle 1/N \rangle\mathrm{.}$$
This isomorphism sends the subgroup $\langle 1/N \rangle$ of
$E_{z'}$ to the image of $E_z[N]$ under the isogeny
$E_z\!\To\!E_z\big/\langle 1/N \rangle$. Also, it sends the basis
$[\,1/p,\,z'/p\,]$ of $E_{z'}[p]$ to the image of the basis
$$[\,(d + c\,z)/p,\,(N^{-1}b + a\,z)/p\,] \, = \,
[\,1/p,\,z/p\,]\,N^{-1}\,\hat{\gamma}$$
of $E_z[p]$. In the cyclotomic case, $d=a$, $c=p$\, and~\,$b$ is a
multiple of $p$, so that $a^2$ equals $N^{-1}$ mod $p$\, \,and\, the 
matrix $\sqrt{N}^{\,-1}\hat\gamma$\, is
trivial in $\PSL_2(\F_p)$. In the non-cyclotomic case, we have 
$\gamma = \pm N\,\hat V$. This completes the proof. \qed

\newpage

\begin{coro}\label{w-rational}
The involution $w$ is defined over $\Q$.
\end{coro}

\noindent\emph{Proof.}\, Take any automorphism $\sigma$ in $\G_\Q$.
Since $w_{\scriptscriptstyle N}$
is defined over $\Q$, $^\sigma\! w$ is still
an involution in $\mathcal{W}(N,p)\setminus \mathcal{G}(N,p)$.
Let $P$ be a non-CM point in $X(N,p)(\Qbar)$ given by a triple
$\left(E,\,C,\,[T_1,T_2]_V\right)$. For a fixed model of the elliptic
curve $E/C$, an isogeny $\lambda\,\colon E\!\To\!E/C$\, with kernel~$C$\, is
determined up to a sign. One has the conjugate isogeny
\,$^\sigma\!\lambda\colon{^\sigma\! E}\!\To\!{^\sigma\!(E/C)}$.
Using Proposition~\ref{w-action} and the isomorphism
$^\sigma\! E/^\sigma\! C\!\To\!{^\sigma\!(E/C)}$ induced by $^\sigma\!\lambda$,
we can verify case by case that ${^\sigma\!}P$ has the same image by
both $w$ and
$^\sigma\! w$. Consider, for instance, the cyclotomic case and assume
\,$^\sigma\!\zeta_p=\zeta_p^{\,r^2}$\, for some $r$ in $\F_p^{\,*}$.
Then, the point ${^\sigma\!}P$ is sent to the isomorphism class of the
triple given by the elliptic curve $^\sigma\!(E/C)$, the cyclic group
$^\sigma\!\lambda(^\sigma\! E[N])$ and the $H_V$--orbit of the basis
$$[\,(r\,\sqrt{N}\,)^{-1}\,\,^\sigma\!\lambda(^\sigma T_1),
\,(r\,\sqrt{N}\,)^{-1}\,\,^\sigma\!\lambda(^\sigma T_2)\,]\mathrm{.}$$
By Proposition~\ref{G(N,p)-action}, this implies
that the matrix in~$\PSL_2(\F_p)$ corresponding to the automorphism
$^\sigma\! w\,w$\, in \,$\mathcal{G}(N,p)$ is the identity, and the result
follows. \qed

\begin{rem}\label{W(N,p)-kp}{\rm \
We can conclude that the Galois covering
\mbox{$X(N,p)\!\To\!X^+(N)$} is defined over $k_p$.
In other words, the function field $k_p\left(X(N,p)\right)$ is a Galois
extension of $k_p\left(X^+(N)\right)$, with group (anti)isomorphic
to~$\mathcal{W}(N,p)$. As a matter of fact, $k_p\left(X(N,p)\right)$
is a Galois extension of $\Q\left(X^+(N)\right)$.
}\end{rem}

Let us finish this section by reviewing the moduli interpretation
of the rational points on $X^+(N)$. The non-cuspidal points of $X_0(N)(\C)$
are in bijection with the isomorphism classes of pairs $(E,C)$, where
$E$\, is a complex elliptic curve and $C$\, is a cyclic subgroup of $E(\C)$
of order $N$\!. Such a point is defined over a number field~$L$ if and
only if $E$ and $C$ are defined over $L$, which means
that \,${^\sigma\!}E=E$\, and \,${^\sigma\!}C=C$\, for all $\sigma$ in $\G_L$.
A point on $X(N,p)$ given by a triple $\left(E,\,C,\,[T_1,T_2]_V\right)$
has image on $X_0(N)$ given by the pair $(E,C)$. In particular, the involution
$w_{\scriptscriptstyle N}$ sends this pair to~$(E/C,\,E[N]/C)$.

Let $E$ be an elliptic curve defined over $L$, and let $E\!\To\!E'$ be an
isogeny with cyclic kernel $C$ of order $N$\!. Assume that~$E$ has no CM,
so that an isogeny from~$E$ to any elliptic curve is determined
up to a sign by its degree. Then, the subgroup~$C$ is defined
over $L$ if and only if $E'$ admits a model over $L$.

Now, suppose that $E$ and $E'$ are defined over a quadratic field $k$,
so that the pair $(E,C)$ defines a $k$-rational point $P$ on $X_0(N)$.
This point is rational if and only if both $E$ and $E'$ have a
model over $\Q$. In this case we say that the couple $\{E,E'\}$ is a
\emph{fake $\Q$-curve of degree $N$}. Otherwise, the image
of~$P$ on~$X^+(N)$ is rational if and only if $E'$ is
isomorphic to the Galois conjugate~\,${^\nu\!}E$\, of~$E$. Indeed, 
since $E$ has no~CM, an isogeny $\mu\colon E\!\To\!{^\nu\!}E$ with 
kernel $C$ sends $E[N]$ to ${^\nu\!}C$, so the existence of such an
isogeny $\mu$\, amounts to the equality \,$w_{\scriptscriptstyle N}(P)={^\nu\!}P$\,
in~\,$X_0(N)(k)$. Thus, every non-cuspidal non-CM rational point on~$X^+(N)$
comes from a pair $(E,C)$ on~$X_0(N)$ defined over some quadratic field
and yielding a (possibly fake) $\Q$-curve of degree~$N$\!.

%---------------------------------------------------------------------
\section{\large The twisted curves in the cyclotomic case}\label{cyclotomic}

Assume $N$ to be a square mod $p$. The structure of this section is
as follows. We first obtain from a modular point of view the fixed field 
of the Galois representation $\varrho_E$ attached in Section~\ref{Weil} to a 
$\Q$-curve $E$ of degree~$N$\!. Next, we produce the twisted modular curves
whose non-cuspidal non-CM rational points give the $\Q$-curves of
degree $N$ realizing a fixed projective mod $p$\, Galois representation
with cyclotomic determinant. We also include a result on the finiteness of
the number of such \mbox{$\Q$-curves}.

Recall that $X^+(N,p)$ denotes the quotient of $X(N,p)$ by the involution $w$. 
The induced map \,$X^+(N,p)\!\To\!X^+(N)$\, is a Galois covering 
with automorphism group $\mathcal{G}(N,p)$ and hence defined over $k_p$ (cf. Proposition \ref{W} 
and Remark~\ref{W(N,p)-kp}). The function field $k_p\left(X^+(N,p)\right)$ is in 
fact a Galois extension of~$\Q\left(X^+(N)\right)$, with group $\PGL_2(\F_p)$\,:\\[-10pt]
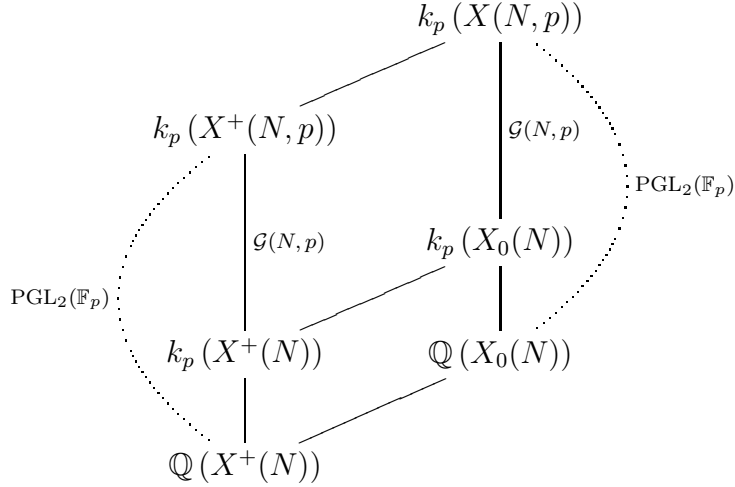
$$\xymatrix{
 & k_p\left(X(N,p)\right) \ar@{-}[dd]^{\mathcal{G}(N,\,p)} \ar@{-}[ld]
 \ar@/^4pc/@{.}[3,0]^{\PGL_2(\F_p)} \\
 k_p\left(X^+(N,p)\right) \ar@{-}[dd]^{\mathcal{G}(N,\,p)}
 \ar@/_4pc/@{.}[3,0]_{\PGL_2(\F_p)} & \\
 & k_p\left(X_0(N)\right) \ar@{-}[d] \ar@{-}[ld] \\
 k_p\left(X^+(N)\right) \ar@{-}[d] & \Q\left(X_0(N)\right) \ar@{-}[ld] \\
 \Q\left(X^+(N)\right) &
}$$
\begin{prop}\label{Weilcyclo}
The function field $k_p\left(X^+(N,p)\right)$ produces, by specialization
over a rational point on $X^+(N)$ corresponding to a \mbox{$\Q$-curve}~$E$, 
the fixed field of the Galois representation $\varrho_E$.
\end{prop}

\noindent\emph{Proof.}\, Let $E$ be a $\Q$-curve of degree $N$ defined
over a quadratic field $k$. Fix an automorphism $\nu$ in $\G_\Q\!\setminus\!\G_k$
and an isogeny $\mu\colon E\!\To\!{^\nu\!}E$ of degree $N$\!. 
%The kernel~$C$ of $\mu$ is defined over $k$, and one has $\mu(E[N])={^\nu\!}C$.
If we let~$C$ be the kernel of $\mu$, the pair~$(E,C)$ 
defines a $k$-rational point on $X_0(N)$ with rational image on $X^+(N)$. 
The preimages on $X^+(N,p)$ of this rational point are
given by the couples $\{P,w(P)\}$ for all points~$P$ on $X(N,p)$ represented
by a triple of the form $\left(E,\,C,\,[T_1,T_2]_V\right)$. If we denote
by~$H$ the subgroup of $\G_{k_p}$ fixing those couples,
what the proposition asserts is that $H$ equals the kernel of~$\varrho_E$. 
This kernel is indeed a subgroup of~$\G_{k_p}$ because
the fixed field of \,$\det\varrho_E$\, is $k_p$ (cf. Corollary \ref{det-proj}).
For a point~$P$ as above, $w(P)$ is given by the triple
$({^\nu\!}E,\,{^\nu\!}C,\,[\sqrt{N}^{\,-1}\mu(T_1),\sqrt{N}^{\,-1}\mu(T_2)]_V)$.
Take now any $\sigma$ in $\G_{k_p}$, so that
$^\sigma\!\zeta_p=\zeta_p^{\,r^2}$ for some $r$ in $\F_p^{\,*}$.
If $\sigma\in\G_k$, then $\sigma\in H$ if and only if
\,${^\sigma\!}P=P$\!, namely if and only if \,${^\sigma T}=\pm r\,T$\,
for all points~$T$ in $E[p]$. If $\sigma\not\in\G_k$, then
$\sigma\in H$ if and only if \,${^\sigma\!}P=w(P)$, namely if and only
if \,${^\sigma T}=\pm r\,\sqrt{N}^{\,-1}\mu(T)$\, for all points $T$
in $E[p]$. Therefore, the result follows from the definition of $\varrho_E$.\qed\\

Suppose that we are now given a Galois representation
$$\varrho\,\colon\,\G_\Q\,\To\,\PGL_2(\F_p)$$
with cyclotomic determinant, which means that the fixed field of \,$\det\varrho$\, is $k_p$. 
For the moduli problem of classifying the $\Q$-curves 
of degree $N$ realizing $\varrho$, we twist the curve $X^+(N,p)$ by certain elements
in the cohomology set $H^1\!\left(\G_\Q, \mathcal{G}(N,p)\right)$. Recall that the twists 
of a curve defined over $\Q$, up to $\Q$-isomorphism, are in bijection with the elements 
in the first cohomology set of~$\G_\Q$ with values in the automorphism group of the curve.

The Galois action on $\mathcal{G}(N,p)$ is known from Corollary~\ref{gal-on-G(N,p)}.
Now, the action by conjugation of \,$\PGL_2(\F_p)$ makes this group 
isomorphic to the automorphism group of $\PSL_2(\F_p)$. Hence, the canonical isomorphism
\,$\mathcal{G}(N,p)\simeq\PSL_2(\F_p)$\, induces an isomorphism
\,$\Aut\left(\mathcal{G}(N,p)\right)\simeq\PGL_2(\F_p)$\, through which the Galois action
on $\mathcal{G}(N,p)$ can be described by the morphism
$$\eta\,\colon\,\G_\Q\,\To\,\Gal(k_p/\Q)\simeq\langle
\hat V\rangle\,\hookrightarrow\,\PGL_2(\F_p)\mathrm{.}$$

Consider then the $1$-cocycles \,$\xi=\varrho_*\eta$\, and 
\,$\xi'=\varrho'_*\eta$, where \,$\varrho_*(\sigma) = {^\mathrm{t}\!}\varrho(\sigma^{-1})$\,
and \,$\varrho'_*(\sigma)=\hat V\varrho_*(\sigma)\,\hat V$\, for all $\sigma$ in $\G_\Q$. The cyclotomic
hypothesis allows us to regard them, through the above canonical isomorphism, as cocycles
with values in~$\mathcal{G}(N,p)$. The cocycle condition for $\xi$, namely
\,$\xi_{\sigma\tau}\,=\,\xi_\sigma\,{^\sigma\!}\xi_\tau$\, for all $\sigma,\tau$ in $\G_\Q$,
can be easily checked case by case, depending on whether \,$\sigma$ and~\,$\tau$ belong to~$\G_{k_p}$ 
or not. The same holds for $\xi'$\!. The cocycle $\xi$ defines a rational model \,$X^+(N,p)_\varrho$\,
for the corresponding twist of $X^+(N,p)$, together with an isomorphism
$$\psi_+\,\colon\,X^+(N,p)_\varrho\,\To\,X^+(N,p)$$
satisfying \,$\psi_+\,=\,\xi_\sigma{^\sigma\!}\psi_+$\, for every~$\sigma$ in~$\G_\Q$. Let us denote 
by \,$X^+(N,p)'_\varrho$\, and~\,$\psi'_+$\, the analogous twist and isomorphism defined by the cocycle $\xi'$\!.

\begin{thm}\label{cyclo1}
There exists a {\rm(}possibly fake{\rm )} $\Q$-curve of degree $N$ realizing $\varrho$ if and only 
if the set of non-cuspidal non-CM rational points on the curves $X^+(N,p)_\varrho$ and 
\,$X^+(N,p)'_\varrho$ is not empty. In this case, the compositions of the isomorphisms \,$\psi_+$ and 
\,$\psi'_+$ with the natural map $X^+(N,p)\!\To\!X^+(N)$ define a surjective map from this set of points 
to the set of isomorphism classes of\,:\\[4pt]
\hspace*{5truemm} $\cdot$ \,$\Q$-curves of degree~$N$ up to Galois conjugation
realizing $\varrho$,\\[4pt]
\hspace*{5truemm} $\cdot$ \,fake $\Q$-curves of degree $N$
realizing~$\varrho$.\\[5pt]
This map is bijective if and only if the centralizer in $\PGL_2(\F_p)$ of the image of $\varrho$ 
is trivial.
\end{thm}

\noindent\emph{Proof.}\, The rational points on $X^+(N,p)_\varrho$\,
correspond via \,$\psi_+$ to the couples of the form $\{P, w(P)\}$,
where $P$ is an algebraic point on $X(N,p)$ such that, for each given automorphism
$\sigma$ in~$\G_\Q$, either \,\,$\xi_\sigma({^\sigma\!}P)=P$\,\, or \,\,$\xi_\sigma({^\sigma\!}P)=w(P)$.

Let $P$ be a non-CM point in $X(N,p)(\Qbar)$ given by a triple
$\left(E,\,C,\,[T_1,T_2]_V\right)$\,. We use the basis $[T_1,T_2]$
of $E[p]$ to fix the isomorphism $\Aut\left(E[p]\right)\simeq\GL_2(\F_p)$.
By virtue of Proposition \ref{G(N,p)-action}, the condition 
\,\mbox{${^\sigma\!}P = \xi_\sigma^{-1}(P)$}\, for all $\sigma$ in $\G_\Q$\, amounts to 
saying that $\{E,E/C\}$ is a fake $\Q$-curve of degree $N$ such that the equality
\begin{equation}\label{xi(N,p)}
\varrho_E(\sigma) \ = \ \left(\begin{matrix}0 & 1\\ 1
& 0\end{matrix}\right) \varrho(\sigma) \left(\begin{matrix}0 & 1\\
1 & 0\end{matrix}\right)
\end{equation}
holds in $\PGL_2(\F_p)$ for every $\sigma$ in $\G_\Q$. Here, we extend the 
notation $\varrho_E$ to the case of elliptic curves over $\Q$ by putting 
$\varrho_E=\overline\rho_E$.

If, on the other hand, there exists \,$\nu$\, in $\G_\Q$ for which
\,$\xi_\nu({^\nu\!}P)=w(P)$, then~$E$ must be a quadratic $\Q$-curve for the point 
$\psi_+^{-1}\left(\{P,w(P)\}\right)$ on $X^+(N,p)_\varrho$ to be rational. Indeed, in 
this case the subgroup of $\G_\Q$ consisting of those automorphisms~$\sigma$\, which
satisfy \,$\xi_\sigma({^\sigma\!}P)=P$\, has index two, so it is of the form $\G_k$ for 
some quadratic field~$k$, and then the condition \,${^\sigma\!}P = \xi_\sigma^{-1}(P)$\, for
all $\sigma$ in $\G_k$ forces the elliptic curve~$E$ and the subgroup~$C$ to be
defined over $k$, while the condition \,$w({^\nu\!}P) = \xi_\nu^{-1}(P)$\, gives
an isogeny $\lambda\,\colon {^\nu\!}E\!\To\!E$ with kernel~${^\nu\!}C$.

So assume now $E$ and $C$ to be defined over a quadratic field $k$\, and let 
$\lambda$ be an isogeny as above. Then, for $\sigma\not\in\G_k$, the 
point $w(^\sigma\! P)$ is represented by the triple given by the elliptic 
curve~$E$, the cyclic group $C$ and the $H_V$--orbit of the basis

\begin{table}[h!]
\centering
\begin{tabular}{ll}
\toprule
$[\,(r\,\sqrt{N}\,)^{-1}\,\lambda({^\sigma}T_1),
\,(r\,\sqrt{N}\,)^{-1}\,\lambda({^\sigma}T_2)\,]$ \ & \
if \ \,${^\sigma\!}\zeta_p = \zeta_p^{\,r^2}$ \\[2pt]
\midrule
$[\,(v\,r\,\sqrt{N}\,)^{-1}\,\lambda({^\sigma}T_1),
\,(v\,r\,\sqrt{N}\,)^{-1}\,\lambda({^\sigma}T_2)\,]\,V$ \ & \
if \ \,${^\sigma\!}\zeta_p = \zeta_p^{\,v\,r^2}$ \\[2pt]
\bottomrule
\end{tabular}
\end{table}
\noindent This comes from Proposition~\ref{w-action} and the
isomorphism \,${^\nu\!}E/{^\nu\!}C\!\To\! E$\, induced by the isogeny $\lambda$. Notice
that the second case does not occur whenever $k=k_p$.

On the other hand, the automorphism $\xi_\sigma^{-1}$
is given by ${^\mathrm{t}\!}\varrho(\sigma)$, if~$\sigma\in\G_{k_p}$, or by
$\hat V\,\,{^\mathrm{t}\!}\varrho(\sigma)$, if $\sigma\not\in\G_{k_p}$.
Then, by applying Proposition \ref{G(N,p)-action} to each
case, we obtain that the point $\psi_+^{-1}\left(\{P,w(P)\}\right)$ on
$X^+(N,p)_\varrho$ is rational if and only if condition~(\ref{xi(N,p)}) holds for 
every $\sigma$ in $\G_\Q$.

Similarly, consider a point on \,$X^+(N,p)'_\varrho$\, corresponding
via \,$\psi'_+$ to a point on~$X^+(N,p)$ obtained from a triple
$\left(E,\,C,\,[T_1,T_2]_V\right)$. By the same reasoning
as above, this point is rational if and only if the pair $(E,C)$
represents a (possibly fake)\, $\Q$-curve of degree $N$ such that,
for every $\sigma$ in $\G_\Q$,
\begin{equation}\label{xi'(N,p)}
\varrho_E(\sigma) \ = \ V \left(\begin{matrix}0 & 1\\
1
& 0\end{matrix}\right) \varrho(\sigma) \left(\begin{matrix}0 & 1\\
1 & 0\end{matrix}\right) V\mathrm{.}
\end{equation}

Let us now consider a (possibly fake) $\Q$-curve of degree $N$ given by some
point~$(E,C)$ on $X_0(N)$ and assume~$\varrho_E = \varrho$. Since this is an 
equality up to conjugation in $\PGL_2(\F_p)$, it amounts to the existence of 
a basis $[T_1,T_2]$ of~$E[p]$ for which condition~(\ref{xi(N,p)}) holds
for every $\sigma$ in $\G_\Q$. Moreover, we can suppose that such a basis is 
sent to either \,$\zeta_p$\, or \,$\zeta_p^{v^{-1}}$ by the Weil pairing. In the 
first case, the image on $X^+(N,p)$ of the triple $\left(E,\,C,\,[T_1,T_2]_V\right)$ 
defines through $\psi_+$ a rational point on $X^+(N,p)_\varrho$. In the second case, 
take $[T_1',T_2']=[T_1,T_2]\,V$. For this new basis, which is sent to \,$\zeta_p$
under the Weil pairing, condition~(\ref{xi'(N,p)}) is satisfied for every~$\sigma$ 
in~$\G_\Q$. So the image on $X^+(N,p)$ of the triple $\left(E,\,C,\,[T_1',T_2']_V\right)$ 
defines through $\psi'_+$ a rational point on~$X^+(N,p)'_\varrho$.

This proves the first part of the statement, including the surjectivity of the map whenever it 
is defined. To discuss its injectivity, consider a point $(E,C)$ on~$X_0(N)$ yielding a 
(possibly fake) $\Q$-curve of degree $N$\!. Suppose that one can take two different rational
points on the twists, corresponding (via $\psi_+$ or~$\psi'_+$) to points on~$X^+(N,p)$ 
obtained from two triples of the form $\left(E,\,C,\,[T_1,T_2]_V\right)$.
Three different cases must be distinguished to complete the proof:
\begin{itemize}
\item[$\cdot$] Both rational points are on $X^+(N,p)_\varrho$ if and only if there is a 
non-trivial element~$\gamma$ in $\PSL_2(\F_p)$, representing a basis change in $E[p]$, 
such that
$$\left(\begin{matrix}0 & 1\\ 1 & 0\end{matrix}\right)
\varrho(\sigma)
\left(\begin{matrix}0 & 1\\ 1 & 0\end{matrix}\right) \ = \
\gamma\,\left(\begin{matrix}0 & 1\\ 1 & 0\end{matrix}\right) \varrho(\sigma)
\left(\begin{matrix}0 & 1\\ 1 & 0\end{matrix}\right) \gamma^{-1}$$
for all $\sigma$ in $\G_\Q$.
This amounts to the existence of a non-trivial element in $\PSL_2(\F_p)$ commuting with all 
the elements in the image of $\varrho$.
\item[$\cdot$] The same characterization is obtained if both points lie on $X^+(N,p)'_\varrho$.
\item[$\cdot$] One of the points is on $X^+(N,p)_\varrho$ and the other on $X^+(N,p)'_\varrho$ 
if and only if there exists $\gamma$ in $\PSL_2(\F_p)$ such that
$$V \left(\begin{matrix}0 & 1\\ 1 & 0\end{matrix}\right)\varrho(\sigma)
\left(\begin{matrix}0 & 1\\ 1 & 0\end{matrix}\right) V \ = \
\gamma\,\left(\begin{matrix}0 & 1\\ 1 & 0\end{matrix}\right) \varrho(\sigma)
\left(\begin{matrix}0 & 1\\ 1 & 0\end{matrix}\right) \gamma^{-1}$$
for every $\sigma$ in $\G_\Q$. This amounts to the existence of an element in
$\PGL_2(\F_p)$ not lying in $\PSL_2(\F_p)$ and commuting with all the elements
in the image of~$\varrho$.\qed
\end{itemize}
%This completes the proof of the statement. \qed\\

%As noticed at the end of Section \ref{Q}, $\Q$-curves of degree~$N$
%are given by rational points on $X^+(N)$. 
Note that the set of points in Theorem~\ref{cyclo1}\, is always finite
whenever the genus of $X^+(N)$ is greater than one. One can assure
this for $N>131$\,: indeed, the modular curve $X_0(N)$ has genus at
least two and is neither hyperelliptic \cite{Ogg} nor bielliptic \cite{Bars}
for any such integer~$N$\!. Using Proposition~\ref{X+(N,p)rational}, one
actually gets the following improvement.

\begin{coro}
\,For \,$N$ square mod $p$, the number of isomorphism classes of\,
\,\mbox{$\Q$-curves} of degree \,$N$ realizing \,$\varrho$\, is finite, except possibly in the case
\mbox{$N=4$, $p=3$.}
\end{coro}

For $\Q$-curves of degree $N$ realizing $\varrho$, a different
moduli description that gets rid of fake $\Q$-curves
can be given for every quadratic field of definition. In order to do
that, we twist $X(N,p)$ by two certain elements in the cohomology set
$H^1\!\left(\G_\Q, \mathcal{W}(N,p)\right)$ that are naturally obtained
from the above cocycles $\xi$ and~$\xi'$ as follows. By Proposition \ref{W} 
and Corollary~\ref{w-rational}, the $\G_\Q$-group $\mathcal{W}(N,p)$
equals the direct product of $\G_\Q$-groups
\,$\mathcal{G}(N,p)\times\langle w\rangle$.
Then, $H^1\!\left(\G_\Q, \mathcal{W}(N,p)\right)$ is also the
direct product of the corresponding cohomology sets.
Fix now a quadratic field~$k$\, and take the Galois character
$$\chi_k\,\colon\,\G_\Q\,\To\,\Gal(k/\Q)\simeq\langle w\rangle\mathrm{.}$$
We then consider the $1$-cocycle \,\,$\xi\,\chi_k$\, and the rational model
\,$X(N,p)_{\varrho,\,k}$\, for the corresponding twist, along with the isomorphism
$$\psi_k\,\colon\,X(N,p)_{\varrho,\,k}\,\To\,X(N,p)$$
satisfying $\psi_k\,=\,(\xi\,\chi_k)_\sigma\,{^\sigma\!}\psi_k$
for every~$\sigma$ in~$\G_\Q$. Analogously, let us denote by
\,$X(N,p)'_{\varrho,\,k}$\, and \,$\psi'_k$\, the twist and the isomorphism
defined by the cocycle~$\xi'\chi_k$\,.

\begin{thm}\label{cyclo2}
There exists a $\Q$-curve of degree $N$ defined over $k$
realizing~$\varrho$ if and only if the set of non-cuspidal non-CM
\,rational points on the curves $X(N,p)_{\varrho,\,k}$\, and
\,$X(N,p)'_{\varrho,\,k}$\, is not empty. In this case, the
compositions of the isomorphisms \,$\psi_k$\, and \,$\psi'_k$\, with the
natural map $X(N,p)\!\To\! X_0(N)$ define a surjective map from this
set of points to the set of isomorphism classes of \,$\Q$-curves
of degree~$N$ defined over~$k$ realizing~$\varrho$. This map is
bijective if and only if the centralizer in $\PGL_2(\F_p)$ of the
image of $\varrho$ is trivial.
\end{thm}

\noindent\emph{Proof.}\, The rational points on $X_{\varrho,\,k}(N,p)$
correspond via~$\psi_k$ to the algebraic points $P$ on $X(N,p)$ such that
$$\xi_\sigma^{-1}(P) \ = \
\left\{\begin{array}{ll}
^\sigma\! P \ & \ \mathrm{for} \ \ \sigma\in\G_k\,\mathrm{,} \\[6pt]
w(^\sigma\! P) \ & \ \mathrm{for} \ \ \sigma\not\in\G_k\,\mathrm{.} \\
\end{array}\right.$$
The proof runs then in a very similar way to that of Theorem \ref{cyclo1},
so we omit the details. In the current case, a non-CM point on 
\,$X(N,p)_{\varrho,\,k}$\, corresponding
via~\,$\psi_k$\, to a triple $\left(E,\,C,\,[T_1,T_2]_V\right)$ is rational if
and only if $E$ is defined over $k$, there exists an isogeny from $E$ to
its Galois conjugate with kernel~$C$ and condition (\ref{xi(N,p)}) holds
for every $\sigma$ in $\G_\Q$ whenever one uses the basis $[T_1,T_2]$ to fix
the isomorphism $\Aut\left(E[p]\right)\simeq\GL_2(\F_p)$. The same characterization 
is valid for the rational points on \,$X(N,p)'_{\varrho,\,k}$\, if we replace \,$\psi_k$\,
by \,$\psi'_k$\, and condition (\ref{xi(N,p)}) by condition (\ref{xi'(N,p)}). \qed

\begin{rem}\label{non-iso1}
{\rm  \ One can check that \,$\xi$\, and \,$\xi'$\, are cohomologous as $1$-cocycles
with values in \,$\mathcal{G}(N,p)$ if and only if the centralizer in $\PGL_2(\F_p)$ of 
the \mbox{image} of~$\varrho$ does not lie in $\PSL_2(\F_p)$. Thus, the twists 
\,$X^+(N,p)_\varrho$\, and \,$X^+(N,p)'_\varrho$\, are not a priori isomorphic over~$\Q$.
The same holds for the twisted curves \,$X(N,p)_{\varrho,\,k}$\, and~\,$X(N,p)'_{\varrho,\,k}$\,.
Moreover, it can be shown that the involution \,$w$\, does not switch the rational points on
\,$X(N,p)_{\varrho,\,k}$\, and \,$X(N,p)'_{\varrho,\,k}$\,, so finding the underlying 
\mbox{$\Q$-curves} requires in general the rational points on both twists.}
\end{rem}

%\begin{rem}\label{non-switch}
%{\rm The action of $w$ does not switch the rational points on the curves 
%$X(N,p)_{\varrho,\,k}$ and $X(N,p)'_{\varrho,\,k}$. Indeed, consider on one of 
%them a rational point corresponding to a triple $\left(E,\,C,\,[T_1,T_2]_V\right)$ 
%on~$X(N,p)$. The image of this triple by~$w$ is
%$({^\nu\!}E,\,{^\nu\!}C,\,[\sqrt{N}^{\,-1}\mu(T_1),\,\sqrt{N}^{\,-1}\mu(T_2)]_V)$,
%where~$\mu$ is an isogeny from~$E$ to its Galois conjugate ${^\nu\!}E$ with kernel~$C$. 
%If we choose the $p$-torsion bases $[T_1,T_2]$ and $[\mu(T_1),\mu(T_2)]$ to fix the 
%isomorphisms $\Aut(E[p])\simeq\GL_2(\F_p)$ and \mbox{$\Aut({^\nu\!}E[p])\simeq\GL_2(\F_p)$,}
%respectively, we then have $\varrho_E(\sigma)=\varrho_{\,{^\nu\!}E}(\sigma)$
%in $\PGL_2(\F_p)$ for all $\sigma$ in~$\G_\Q$. Hence, according to the proof of
%Theorem~\ref{cyclo2}, the new triple corresponds to a rational point on the same twisted 
%curve as the given point.}
%\end{rem}

%----------------------------------------------------------------------
\section{\large The twisted curve in the non-cyclotomic case}\label{non-cyclotomic}

Assume $N$ to be a non-square mod $p$. This section is the analogue of the 
previous one for the non-cyclotomic case. Unlike the cyclotomic case, now the 
quadratic field of definition for the potential $\Q$-curves of degree $N$
realizing a given projective mod~$p$\, Galois representation is fixed by the 
determinant. Moreover, only one twist is needed for the moduli classification 
of such $\Q$-curves. We prove this in Theorem~\ref{non-cyclo} below.
For the sake of completeness, let us begin as before with the modular construction 
of the fixed field of the Galois representation $\varrho_E$ attached to a $\Q$-curve $E$
of degree $N$\!. The procedure is now more intricate.

Recall that the group $\mathcal{W}(N,p)$ of the covering
$X(N,p)\!\To\! X^+(N)$ is canonically isomorphic to $\PGL_2(\F_p)$. The action 
by conjugation of this group makes it isomorphic to its automorphism group.
Thus, by virtue of Corollary~\ref{gal-on-G(N,p)} and Corollary~\ref{w-rational}, 
the Galois action on~$\mathcal{W}(N,p)$ is given by the morphism
$$\eta\,\colon\,\G_\Q\,\To\,\Gal(k_p/\Q)\simeq\langle
w\rangle\,\hookrightarrow\,\mathcal{W}(N,p)\mathrm{,}$$
where we identify $\mathcal{W}(N,p)$ with its (inner) automorphism group.

Let $\widetilde X(N,p)$ be the twist of $X(N,p)$ defined by the $1$-cocycle $\eta$.
Likewise, denote by $\widetilde X_0(N)$ the twist of $X_0(N)$ defined by the $1$-cocycle
$$\G_\Q\,\To\,\Gal(k_p/\Q)\simeq\langle w_{\scriptscriptstyle N}\rangle\mathrm{.}$$
We write $\widetilde X^+(N)$ for the quotient of $\widetilde X_0(N)$ by
the involution corresponding to~$w_{\scriptscriptstyle N}$. Consider the
following commutative diagram, where the morphisms are the natural ones:
$$
\xymatrix{
 \widetilde X(N,\,p) \ar@{->}[r]^{\simeq} \ar@{->}[d] & X(N,\,p)
\ar@{->}[d]  \\
\widetilde X_0(N) \ar@{->}[r]^{\simeq} \ar@{->}[d] & X_0(N) \ar@{->}[d] \\
\widetilde X^+(N) \ar@{->}[r]^{\simeq} & X^+(N)
}$$
As remarked in the proof of the next lemma, the isomorphism
$\widetilde X^+(N)\!\To\! X^+(N)$ is actually defined over $\Q$.

\begin{lema}\label{Qcovering}
The Galois covering \,$\widetilde X(N,p)\To\widetilde X^+(N)$\, is
defined over $\Q$.
\end{lema}
\noindent\emph{Proof.}\, Denote by
$\phi\,\colon \widetilde X(N,p)\!\To\! X(N,p)$\, and
\,$\phi_0\colon \widetilde X_0(N)\!\To\! X_0(N)$ the isomorphisms in the
above diagram. They are defined over \,$k_p$\, and satisfy
${^\sigma\!}\phi\,\phi^{-1}=w$\, and
\,${^\sigma\!}\phi_0\,\phi_0^{-1}=w_{\scriptscriptstyle N}$\,
for $\sigma\not\in\G_{k_p}$. Then, the involution
$\phi_0^{-1}\,w_{\scriptscriptstyle N}\,\phi_0$ on $\widetilde
X_0(N)$ is defined over $\Q$. Hence, so is the corresponding quotient map
\mbox{$\widetilde X_0(N)\!\To\!\widetilde X^+(N)$.}
The isomorphism $\phi_+\colon X^+(N)\!\To\!\widetilde X^+(N)$ induced by
$\phi_0^{-1}$ sends a couple $\{P, w_{\scriptscriptstyle N}(P)\}$
to~$\{\phi_0^{-1}(P), \phi_0^{-1}\,w_{\scriptscriptstyle N}(P)\}$. It
is easily checked to satisfy ${^\sigma\!}\phi_+=\phi_+$ for all
$\sigma$ in~$\G_\Q$. The same is true for the morphism
%since $w$ extends $w_{\scriptscriptstyle N}$,
\mbox{$\widetilde X(N,p)\!\To\!\widetilde X_0(N)$} induced from the natural map
$X(N,p)\!\To\! X_0(N)$ by the isomorphisms $\phi$\, and~$\phi_0$. Finally,
the automorphisms 
%in the group \,$\phi^{-1}\,\mathcal{W}(N,p)\,\phi$\, 
of the covering \mbox{$\widetilde X(N,p)\!\To\!\widetilde X^+(N)$} are also defined over $\Q$.
%given an automorphism $\vartheta$ in $\mathcal{W}(N,p)$,
Indeed, the relation \,${^\sigma\!}(\phi^{-1}\,\vartheta\,\,\phi) =
\phi^{-1}\,w\,(w\,\vartheta\,w)\,w\,\phi = \phi^{-1}\,\vartheta\,\,\phi$\,\,
holds for~\mbox{$\vartheta\in\mathcal{W}(N,p)$} and~$\sigma\not\in\G_{k_p}$. For another
proof of the existence of such a rational covering we refer to
\cite{Sh74}. \qed
%, where the problem of finding a minimal field of definition for
%$X(N,p)\to X_0(N)$ is considered. \qed

\begin{rem}\label{Weilnoncyclo-ffield}{\rm \ The function field of
$\widetilde X(N,p)$ over $\Q$ is identified, through the isomorphism $\phi$\,
in the proof of Lemma \ref{Qcovering}, with a subfield of $k_p\left(X(N,p)\right)$.
%with the subfield of functions $f$ in $k_p\left(X(N,p)\right)$ 
%satisfying \,${^\sigma\!}f=fw$\, for some (hence all) automorphism $\sigma\not\in\G_{k_p}$. 
%Notice, in particular, that $\Q\,(\widetilde X(N,p))$contains the function field 
%over~$\Q$ of the quotient curve $X(N,p)/w$, 
As shown in the following diagram, it is a Galois extension of \,$\Q\left(X^+(N)\right)$
with group isomorphic to $\PGL_2(\F_p)$\,:%\\[-10pt]
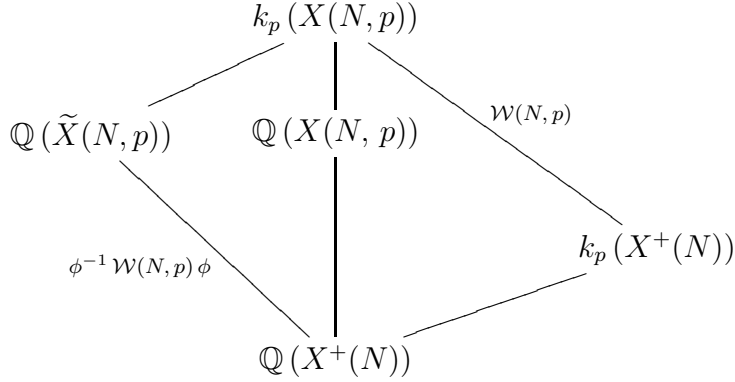
$$\xymatrix{
 & k_p\left(X(N,p)\right) \ar@{-}[ld] \ar@{-}[d]
 \ar@{-}[rrdd]^{\mathcal{W}(N,\,p)} & & \\
 \Q\,(\widetilde X(N,p)) %\ar@{-}[rd]
 %\ar@/_2pc/@{.}[3,1]_{\phi^{-1}\,\mathcal{W}(N,\,p)\,\phi} &
 \ar@{-}[rdd]_{\phi^{-1}\,\mathcal{W}(N,\,p)\,\phi} &
 \Q\left(X(N,\,p)\right) \ar@{-}[dd] & & \\
 %& \Q\left(X(N,p)/w\right) \ar@{-}[dd] & & \\
 &  &  & k_p\left(X^+(N)\right) \ar@{-}[lld] \\
 & \Q\left(X^+(N)\right) & &
}$$
}
\end{rem}

\begin{prop}\label{Weilnoncyclo}
The function field \,$\Q\,(\widetilde X(N,p))$ gives, by specialization
over a rational point on $X^+(N)$ corresponding to a
$\Q$-curve $E$, the fixed field of the representation $\varrho_E$.
\end{prop}

\noindent\emph{Proof.}\, With the same notations as in the proof of
Proposition \ref{Weilcyclo}, take a cyclic isogeny \,$\mu\,\colon E\!\To\!{^\nu\!}E$\,
with kernel $C$ of order $N$\!. Consider the isomorphism 
$\phi\,\colon \widetilde X(N,p)\!\To\! X(N,p)$\, in the proof of Lemma \ref{Qcovering}. 
Let $H$ be the subgroup of $\G_\Q$ fixing the points on $\widetilde X(N,p)$ 
corresponding through $\phi$\, to the points on $X(N,p)$ of the form $P$ or $w(P)$, 
where $P$ is given by a triple of the form $\left(E,\,C,\,[T_1,T_2]_V\right)$. 
We must show that $H$ is the kernel of $\varrho_E$. For a point~$P$ as above, 
$w(P)$ is represented by the triple given by \,${^\nu\!}E$, ${^\nu\!}C$ and the 
$H_V$--orbit of the basis $[\mu(T_1),\mu(T_2)]V$\!.
Using the definition of $\phi$, we see that the group~$H$ consists
of those \,$\sigma\in\G_{k_p}$\, satisfying \,${^\sigma\! P}=P$ and
those \,$\sigma\not\in\G_{k_p}$\, satisfying
\,${^\sigma\! P}=w(P)$. Moreover, any such~$\sigma$ lies in $\G_k$ if and only
if it lies in $\G_{k_p}$. Take now any automorphism $\sigma$ in $\G_\Q$.
If~\,$^\sigma\!\zeta_p=\zeta_p^{\,r^2}$\, for some $r$ in~$\F_p^{\,*}$,
then $\sigma\in H$ if and only if \,${^\sigma\!}P=P$\!, namely if and
only if \,${^\sigma T}=\pm r\,T$ for all points $T$ in $E[p]$.
If~\,$^\sigma\!\zeta_p=\zeta_p^{\,r^2N^{-1}}$ for some $r$ in~$\F_p^{\,*}$, then
$\sigma\in H$ if and only if \,${^\sigma\!}P=w(P)$, namely if and only
if \,${^\sigma T}=\pm r\,N^{-1}\,\mu(T)$ for all points $T$
in $E[p]$. So the result follows from the definition of $\varrho_E$.
\qed\\

%\vspace{10truemm}

Suppose that we have now a Galois representation
$$\varrho\,\colon\,\G_\Q\,\To\,\PGL_2(\F_p)$$
with non-cyclotomic determinant. Recall that any $\Q$-curves of degree $N$
realizing~$\varrho$ must be defined over the fixed field of 
\,$\eps\det\varrho$, where $\eps$ is the character attached to $k_p$ 
(cf. Corollary \ref{det-proj}).
Denote this quadratic field by $k$. For the moduli classification of such
\mbox{$\Q$-curves,} we produce a twist of $X(N,p)$ from a certain element
in the cohomology set $H^1\!\left(\G_\Q, \mathcal{W}(N,p)\right)$, as follows.
The canonical isomorphism \,$\mathcal{W}(N,p)\simeq\PGL_2(\F_p)$\, allows us
to regard the projective representation \,$\varrho_*$\, in Section \ref{cyclotomic}\,
as a morphism taking values in~$\mathcal{W}(N,p)$. As before, let~$\eta$ stand for
the morphism giving the Galois action on~$\mathcal{W}(N,p)$. Then, consider the $1$-cocycle 
%$\G_\Q$ with values in $\mathcal{W}(N,p)$ given by the product
\,$\xi=\varrho_*\eta$. For the twist of $X(N,p)$ defined by \,$\xi$, we fix
a rational model $X(N,p)_\varrho$\, along with an isomorphism
$$\psi\,\colon\,X(N,p)_\varrho\,\To\,X(N,p)$$
satisfying \,$\psi\,=\,\xi_\sigma{^\sigma\!}\psi$\, for every~$\sigma$ in~$\G_\Q$.

\begin{thm}\label{non-cyclo}
There exists a \,$\Q$-curve of degree $N$ realizing $\varrho$ if and only
if the set of non-cuspidal non-CM \,rational points on the curve
$X(N,p)_\varrho$ is not empty. In this case, the composition of \,$\psi$ with
the natural map $X(N,p)\!\To\! X^+(N)$ defines a
surjective map from this set of points to the set of isomorphism classes
of \,\mbox{$\Q$-curves} of degree~$N$ up to Galois conjugation realizing~$\varrho$.
This map is bijective if and only if the centralizer in $\PGL_2(\F_p)$ of
the image of $\varrho$ is trivial.
\end{thm}

\noindent\emph{Proof.}\, The first part of the proof goes along the
lines of those of Theorem~\ref{cyclo1} and Theorem~\ref{cyclo2}. Let us fix an
automorphism \,$\nu$\, in $\G_\Q\!\setminus\!\G_k$.
The rational points on $X_\varrho(N,p)$ correspond via \,$\psi$
to the algebraic points~$P$ on $X(N,p)$ satisfying
$$\varrho_*(\sigma)^{-1}(P) \ = \
\left\{\begin{array}{ll}
^\sigma\! P \ & \ \mathrm{for} \ \sigma\in\G_{k_p}\,\mathrm{,} \\[6pt]
w(^\sigma\! P) \ & \ \mathrm{for} \ \sigma\not\in\G_{k_p}\,\mathrm{.} \\
\end{array}\right.$$
Note that the automorphism $\varrho_*(\sigma)^{-1}$ belongs
to $\mathcal{G}(N,p)$
if and only if $\sigma$ lies in either both $\G_k$ and $\G_{k_p}$ or
none of them. In particular, a non-CM point $P$ given by a triple
$\left(E,\,C,\,[T_1,T_2]_V\right)$ can satisfy the above condition
only if $E$ and $C$ are defined over $k$ and there is an isogeny
$\lambda\,\colon{^\nu\!}E\!\To\! E$ with kernel ${^\nu\!}C$.
With these hypotheses on $E$ and $C$, and for $\sigma\not\in\G_k$,
the point~$w(^\sigma\! P)$ is represented by the triple given by $E$,
\,$C$ and the $H_V$--orbit of the basis
%the choice depends as usual on the cyclotomic behaviour of $\sigma$.
\begin{table}[h!]
\centering
\begin{tabular}{ll}
\toprule
$[\,r^{-1}\,\lambda({^\sigma}T_1),\,r^{-1}\,\lambda({^\sigma}T_2)\,]\,V$ \ & \
if \ \,${^\sigma\!}\zeta_p = \zeta_p^{\,r^2}$ \\[2pt]
\midrule
$[\,r^{-1}\,\lambda({^\sigma}T_1),\,r^{-1}\,\lambda({^\sigma}T_2)\,]$ \ & \
if \ \,${^\sigma\!}\zeta_p = \zeta_p^{\,r^2 N^{-1}}$ \\[2pt]
\bottomrule
\end{tabular}
\end{table}

\noindent In the second case, and also for
$\sigma\in\G_k\cap\G_{k_p}$,
the automorphism~$\varrho_*(\sigma)^{-1}$ is given by the matrix
${^\mathrm{t}\!}\varrho(\sigma)$ in $\PSL_2(\F_p)$. In the other case,
and also for $\sigma\in\G_k\!\setminus\!\G_{k_p}$, the automorphism
$w\,\varrho_*(\sigma)^{-1}$ is given by the matrix
$\hat V\,\,{^\mathrm{t}\!}\varrho(\sigma)$ in $\PSL_2(\F_p)$.
So, taking as in the proof of Theorem~\ref{cyclo1} the basis $[T_1,T_2]$
to fix the isomorphism $\Aut\left(E[p]\right)\simeq\GL_2(\F_p)$, 
condition~(\ref{xi(N,p)}) is again seen to characterize the rationality of 
the point $\psi^{-1}(P)$.

Consider now a \mbox{non-CM} elliptic
curve $E$ defined over $k$\, and an isogeny $\mu\,\colon E\!\To\!{^\nu\!}E$\, with
kernel $C$, and assume~\mbox{$\varrho_E = \varrho$.}
This equality amounts to the existence of a basis $[T_1,T_2]$ of~$E[p]$
for which condition~(\ref{xi(N,p)}) holds for every $\sigma$
in~$\G_\Q$. We can further suppose that such a basis is sent to
either \,$\zeta_p$\, or \,$\zeta_p^{N^{-1}}$ by the Weil pairing. In the first
case, the point~$P$ on~$X(N,p)$ given
by the triple $\left(E,\,C,\,[T_1,T_2]_V\right)$ defines through $\psi$ a
rational point on $X(N,p)_\varrho$. In the second case, the triple
$\left({^\nu\!}E,\,{^\nu\!}C,\,[\mu(T_1),\mu(T_2)]_V\right)$ represents a point on
$X(N,p)$ lying above the same point on~$X^+(N)$ as $P$ and
corresponding via $\psi$ to a rational point on $X(N,p)_\varrho$.
Indeed, if we choose the basis $[\mu(T_1),\mu(T_2)]$
to fix the isomorphism $\Aut({^\nu\!}E[p])\simeq\GL_2(\F_p)$, we get the equality
\,$\varrho_{\,{^\nu\!}E}(\sigma)=\varrho_E(\sigma)$\, for all $\sigma$ in~$\G_\Q$.
%(cf.~Remark~\ref{non-switch}).

Lastly, let us consider two different rational
points on \,$X(N,p)_\varrho$\, corresponding via \,$\psi$\, to non-CM points
$P$ and $Q$ on $X(N,p)$ with the same image on~$X^+(N)$. Let the triple 
$\left(E,\,C,\,[T_1,T_2]_V\right)$ represent the point $P$
and fix an isogeny $\mu\,\colon E\!\To\!{^\nu\!}E$\, with kernel $C$.
%Then, a triple representing $w(P)$ is given by the elliptic curve ${^\nu\!}E$,
%the subgroup ${^\nu\!}C$ and the $H_V$--orbit of the basis
%$[\mu(T_1),\mu(T_2)]\,V$
We must then distinguish two cases for the point $Q$\,: 
\begin{itemize}
\item[$\cdot$] It lies over the pair $(E,C)$ on~$X_0(N)$ if and only if there
is a non-trivial element~$\gamma$ in $\PSL_2(\F_p)$,
representing a basis change in $E[p]$, such that
$$\left(\begin{matrix}0 & 1\\ 1 & 0\end{matrix}\right)
\varrho(\sigma)
\left(\begin{matrix}0 & 1\\ 1 & 0\end{matrix}\right) \ = \
\gamma\,\left(\begin{matrix}0 & 1\\ 1 & 0\end{matrix}\right) \varrho(\sigma)
\left(\begin{matrix}0 & 1\\ 1 & 0\end{matrix}\right) \gamma^{-1}$$
for all $\sigma$ in $\G_\Q$. This amounts to the existence of a
non-trivial element in $\PSL_2(\F_p)$ commuting with all the elements in
the image of $\varrho$.
\item[$\cdot$] Otherwise, a triple representing $Q$ is given by the elliptic curve ${^\nu\!}E$,
the subgroup \,${^\nu\!}C$ and the $H_V$--orbit of a basis obtained from \,$[\mu(T_1),\mu(T_2)]\,V$\,
by a basis change preserving the Weil pairing. Thus, this case amounts to the existence
of an element $\gamma$ in $\PSL_2(\F_p)$
%, representing a basis change in ${^\nu\!}E[p]$, 
such that
$$V \left(\begin{matrix}0 & 1\\ 1 & 0\end{matrix}\right)
\varrho(\sigma)
\left(\begin{matrix}0 & 1\\ 1 & 0\end{matrix}\right) V \ = \
\gamma\,\left(\begin{matrix}0 & 1\\ 1
& 0\end{matrix}\right) \varrho(\sigma) \left(\begin{matrix}0 & 1\\
1 & 0\end{matrix}\right)\,\gamma^{-1}$$
for all $\sigma$ in $\G_\Q$. This is in turn equivalent to the
existence of an element in $\PGL_2(\F_p)\!\setminus\PSL_2(\F_p)$
commuting with all the elements in the image of~$\varrho$. \qed\\
\end{itemize}
 
%We conclude with the following finiteness result
By the same reasoning as in Section \ref{cyclotomic}, the set of points in 
Theorem \ref{non-cyclo} is always finite whenever $N>131$. A stronger result 
is obtained from \mbox{Corollary \ref{X(N,p)rational-elliptic}.}

\begin{coro}
For $N$ non-square mod $p$, the number of isomorphism classes of \,$\Q$-curves 
of degree~$N$ realizing \,$\varrho$ is finite, unless $N=2$ and $p=3$.
\end{coro}

%\begin{rem}\ref{non-CM}????????
%{\rm Assume that $\varrho$ is surjective. If $E$ is a quadratic
%$\Q$-curve realizing $\varrho$, then the image of the projective Galois
%representation $\overline\rho_{E,\,p}$ is not surjective and hence $E$
%cannot have CM. Thus, Theorem \ref{non-cyclo} can be
%restated by saying that the non-cuspidal rational points on the curve
%$X(N,p)_\varrho$ are in bijection with the (isomorphism classes of) pairs
%of $\Q$-curves of degree $N$ realizing $\varrho$. Analogous remarks can be
%made for Theorem \ref{cyclo1} and Theorem \ref{cyclo2}.
%}
%\end{rem}

%--------------------------------------------------------------------------------

\vspace{1truemm}

\section*{\large Acknowledgements}
I wish to express my gratitude to Joan-Carles Lario for his help and encouragement 
through this work. I am most indebted to René Schoof and the \emph{Dipartimento di Matematica}
of the \emph{Universit\`a di Roma "Tor Vergata"} as well as to the 
\emph{D\'epartement de Math\'ematiques de Besançon},
where this research was carried out with financial support from the RTN
European Network \emph{Galois Theory and Explicit Methods in Arithmetic}. 
It is also a pleasure to thank Gabriel Cardona for some helpful comments 
on an earlier version of the paper.

%-------------------------------------------------------------------------------

\vspace{1truemm}

\bibliographystyle{alpha}
\bibliography{imrn}

\vfill

\begin{tabular}{l}
Julio Fern\'andez \\[2pt]
Departament de Matem\`atica Aplicada 4\\[2pt]
Universitat Polit\`ecnica de Catalunya\\[2pt]
EPSEVG, av. V\'ictor Balaguer\\[2pt]
E-08800 Vilanova i la Geltr\'u (Barcelona)\\[2pt]
\texttt{julio@mat.upc.edu}
\end{tabular}

\end{document}